\documentclass[12pt,a4paper]{article}
\usepackage{amsmath, amssymb}
\usepackage{graphicx, subcaption}
\usepackage{geometry}
\usepackage{microtype} 
\usepackage{caption} 
\usepackage{amsthm} 
\usepackage{hyperref} 
\usepackage{algorithm}
\usepackage{algorithmic}
\usepackage{booktabs}
\usepackage{float}

\newtheorem{conjecture}{Conjecture}
\theoremstyle{remark} 

\title{\textbf{A Gas-Driven Algorithm for Variants of the Moving Sofa Problem}}
\author{Xingyi He\thanks{
College of Architecture and Environment, Sichuan University,
Chengdu 610000, China.
E-mail: \texttt{xingyihe@stu.scu.edu.cn}.
}}
\date{\today}

\begin{document}
\maketitle

\begin{abstract}
This paper presents a numerical algorithm based on the idea of dynamical systems for solving the moving sofa problem. 
By introducing a physical model driven by gas pressure, 
we transform the geometric optimization problem into a dynamical system. 
Numerical experiments show that the method effectively approximates the known Gerver's sofa result, 
and provides numerical estimates of the possible maximum area for corridors with different angles. 
In particular, we observe an intersection of the two motion patterns at $43.327\ldots^\circ$, 
where the locally maximal sofa area is $1.8674\ldots$, 
and the dominant pattern switches on either side of this critical angle.
\end{abstract}

\section{Introduction}

The moving sofa problem is a classic optimization problem, first posed by Moser in 1966 \cite{Moser1966}. The problem can be stated as:

\begin{quote}
\textit{What is the planar shape of maximal area that can be moved around a right-angled corner in a corridor of unit width?}
\end{quote}

Although the problem is intuitive to state, the simultaneous translational and rotational motions of a planar shape within a confined space render its quantitative analysis and solution extremely complex. From an optimization perspective, the problem can be viewed as an infinite-dimensional non-convex optimization problem, where the optimization variables are connected planar regions satisfying the geometric constraints of the corridor.

Since its inception, the problem has attracted extensive research. In 1968, Hammersley proposed a shape consisting solely of circular arcs and straight line segments, with an area of $\frac{\pi}{2} + \frac{2}{\pi} \approx 2.2074$, and also gave an upper bound of $2\sqrt{2}$ \cite{Hammersley1968}. Subsequently, numerical methods were employed to explore possible optimal shapes. In 1973, Maruyama designed a program to find approximate solutions to the general sofa problem \cite{Maruyama1973}, and in 1976, Wagner used Monte Carlo algorithms to optimize the right-angled sofa \cite{Wagner1976}.

\begin{figure}[htbp]
\centering
\includegraphics[width=0.8\textwidth]{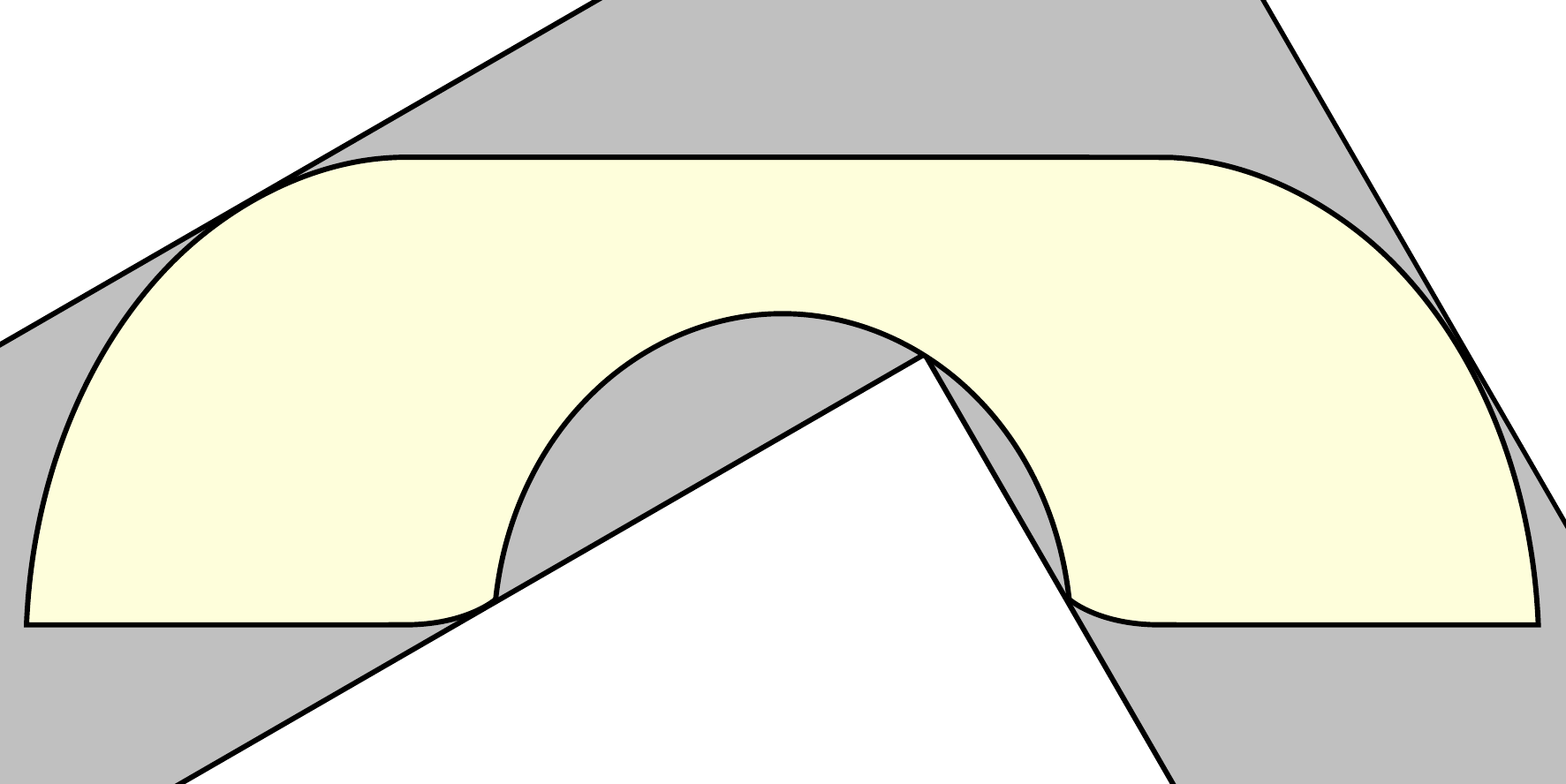}
\caption{Gerver's sofa}
\label{fig:90r}
\end{figure}

In 1992, Gerver derived necessary conditions for a sofa to have maximal area, known as the balanced polygon theory. Based on this, he constructed a more complex sofa shape composed of three straight line segments and fifteen curved segments (see Figure~\ref{fig:90r}), with an area of approximately 2.2195\ldots, which is currently considered the largest known area for a sofa \cite{Gerver1992,Baek2024}. In 2014, Gibbs used numerical methods to solve the right-angled sofa and the double-sided right-angled sofa \cite{Gibbs2014}. In 2016, Romik extended Gerver's approach by transforming the necessary conditions for maximal area into a family of six ordinary differential equations, and using this method he obtained the shape of the double-sided right-angled sofa, which can be expressed in closed form \cite{Romik2016}. In 2018, Kallus and Romik used computer-assisted methods to prove an upper bound of 2.37 for the maximal sofa area, and also showed that a sofa of maximal area must rotate through an angle of at least $81.203^\circ$ when navigating the corner \cite{KallusRomik2018}. In 2024, Deng applied the calculus of variations to the moving sofa problem, while Leng et al. employed deep learning methods for numerical solutions \cite{Deng2024,Leng2024}. Most recently, Baek presented a proof of the optimality of the Gerver's sofa \cite{Baek2024}. Furthermore, Georgiev et al. used AlphaEvolve to numerically solve a three-dimensional variant of the moving sofa problem \cite{Georgiev2025}.

\begin{figure}[htbp]
\centering

\begin{subfigure}{0.32\textwidth}
\centering
\includegraphics[width=\linewidth]{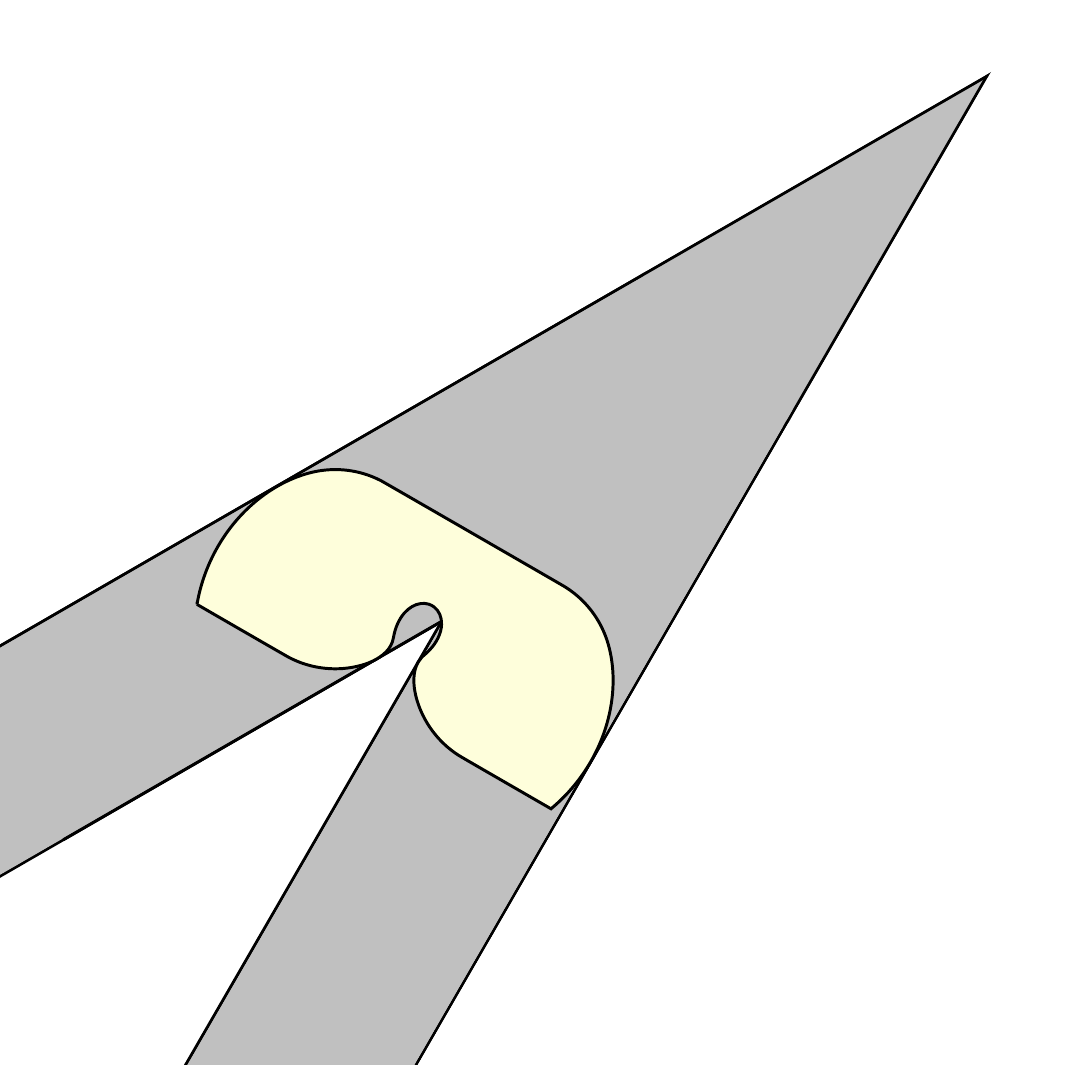}
\end{subfigure}
\begin{subfigure}{0.32\textwidth}
\centering
\includegraphics[width=\linewidth]{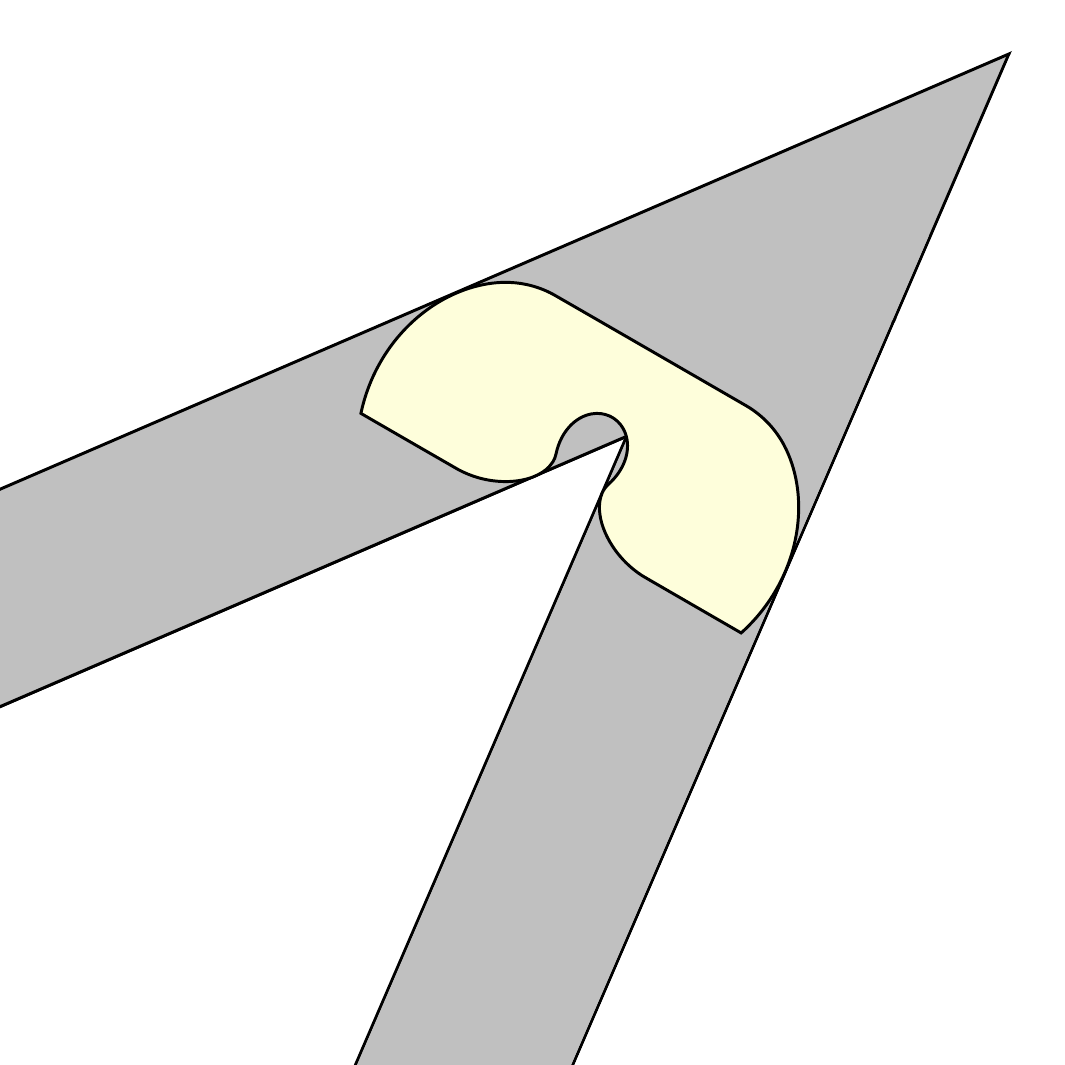}
\end{subfigure}
\begin{subfigure}{0.32\textwidth}
\centering
\includegraphics[width=\linewidth]{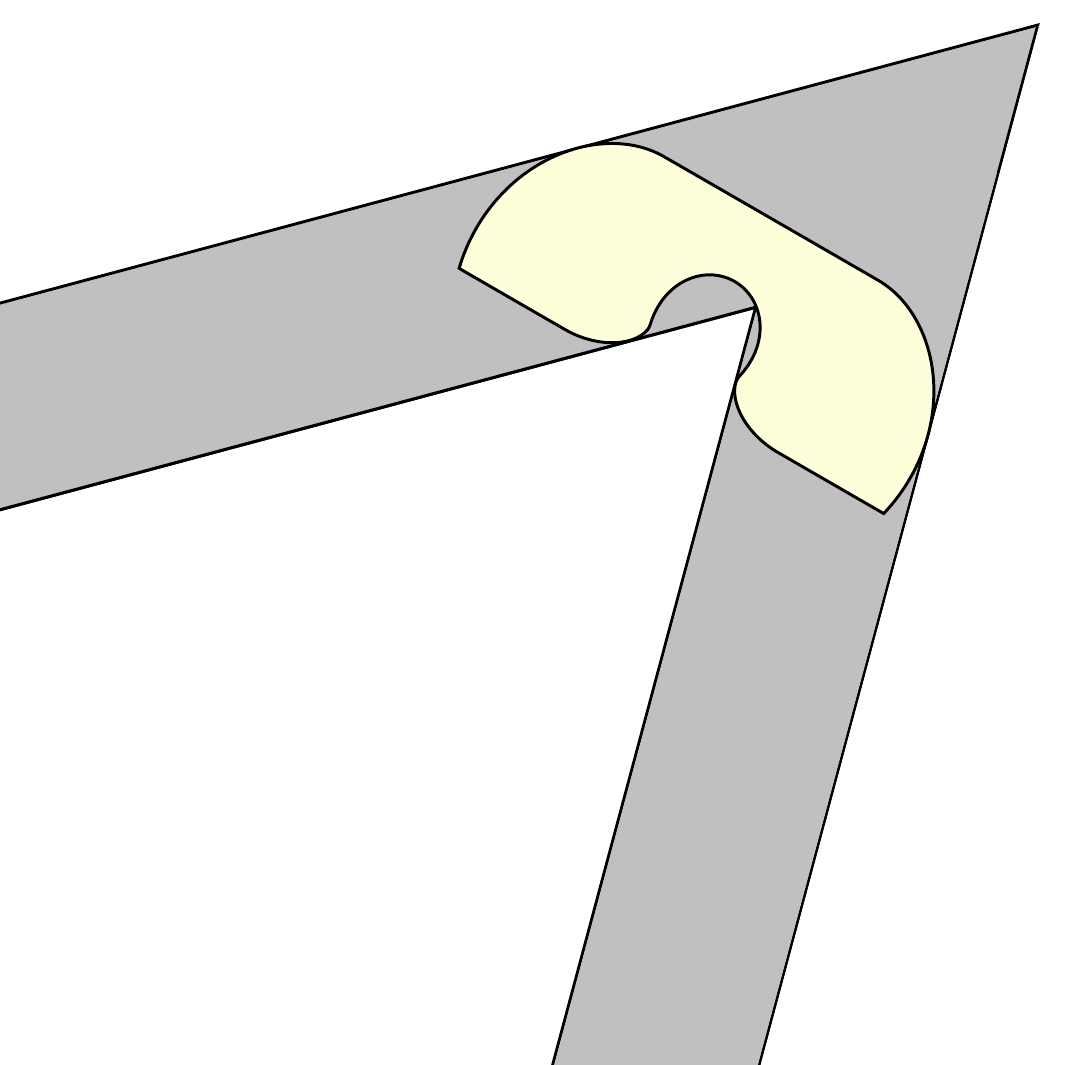}
\end{subfigure}

\vspace{0.2cm}

\begin{subfigure}{0.32\textwidth}
\centering
\includegraphics[width=\linewidth]{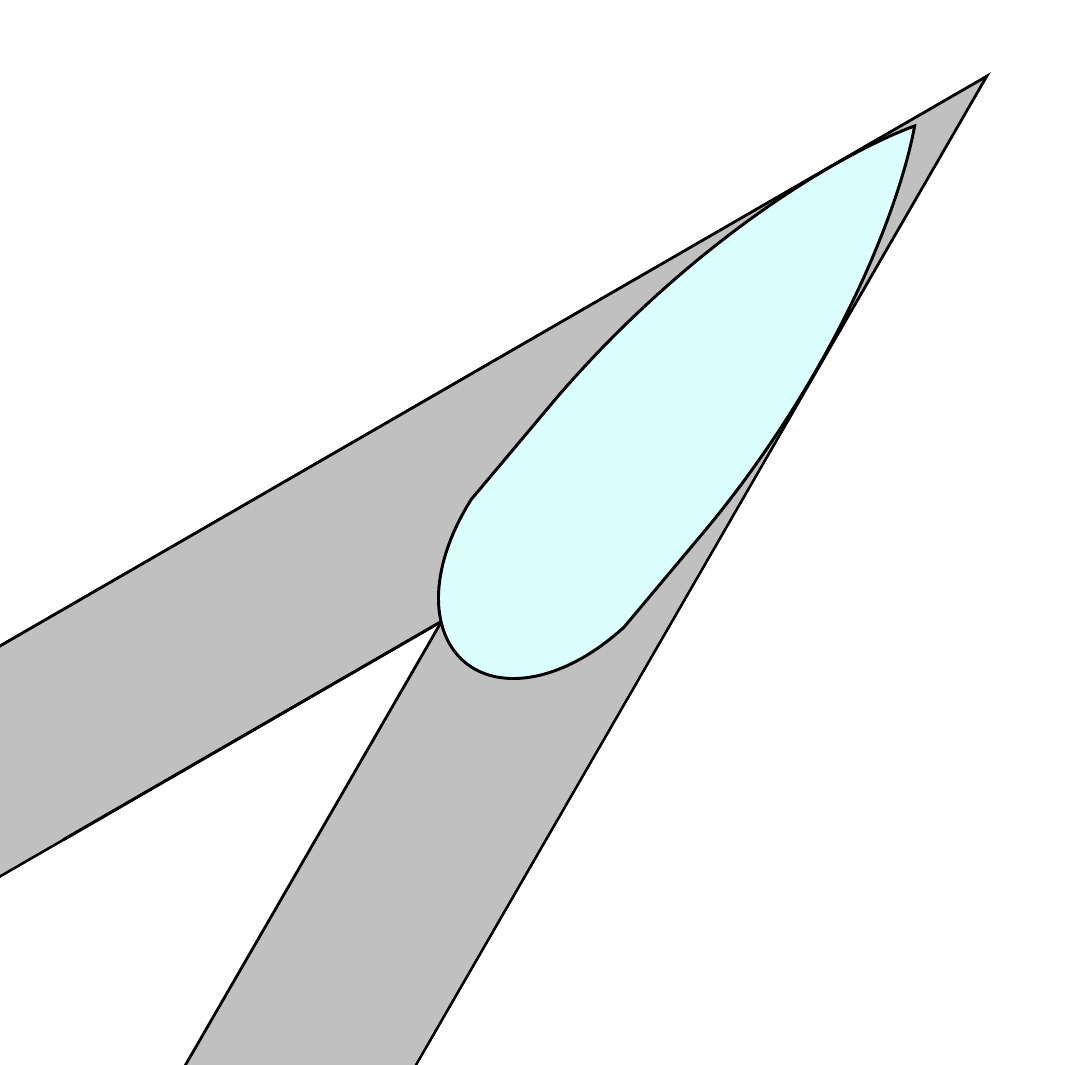}
\caption{$30^\circ$}
\end{subfigure}
\begin{subfigure}{0.32\textwidth}
\centering
\includegraphics[width=\linewidth]{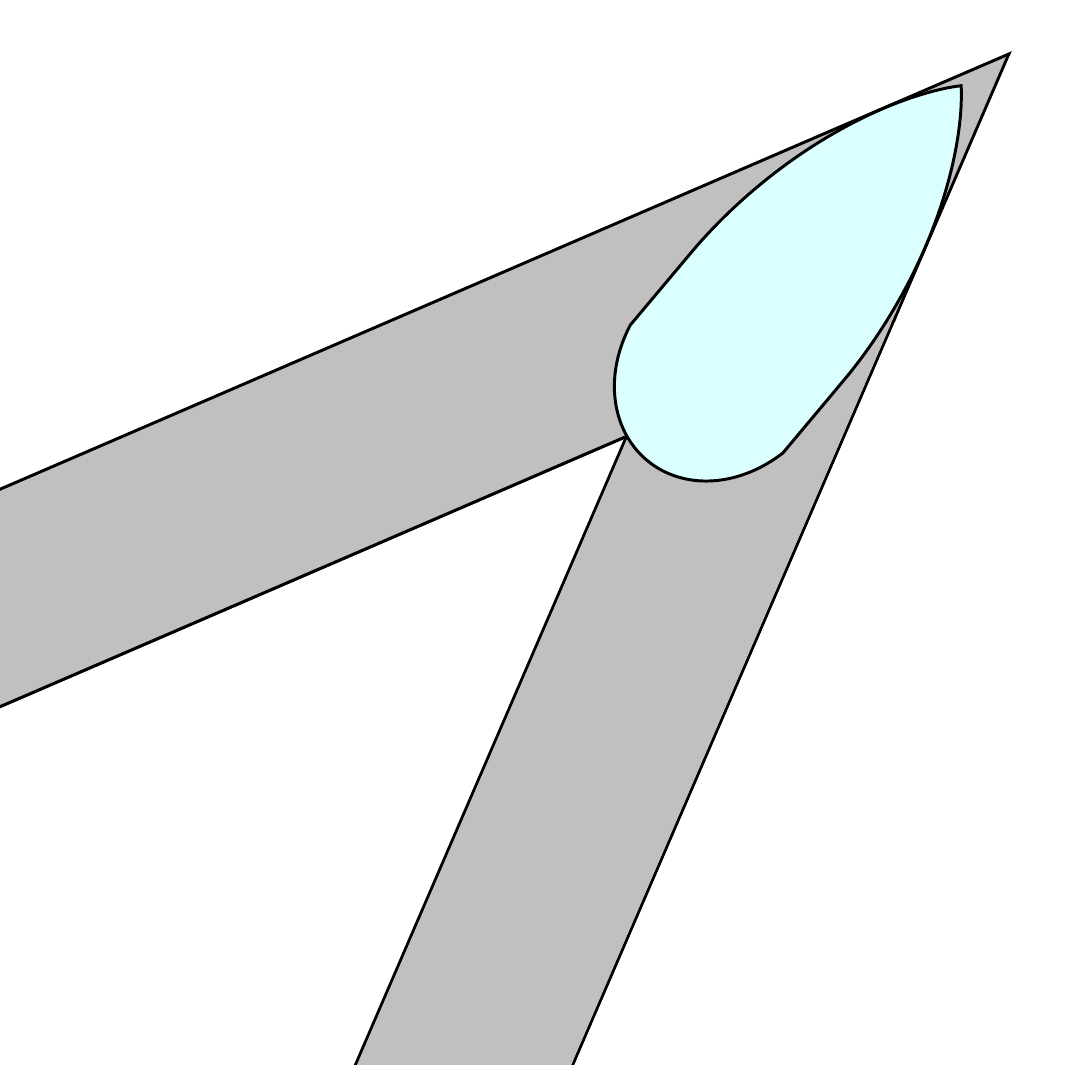}
\caption{$43.327^\circ$}
\end{subfigure}
\begin{subfigure}{0.32\textwidth}
\centering
\includegraphics[width=\linewidth]{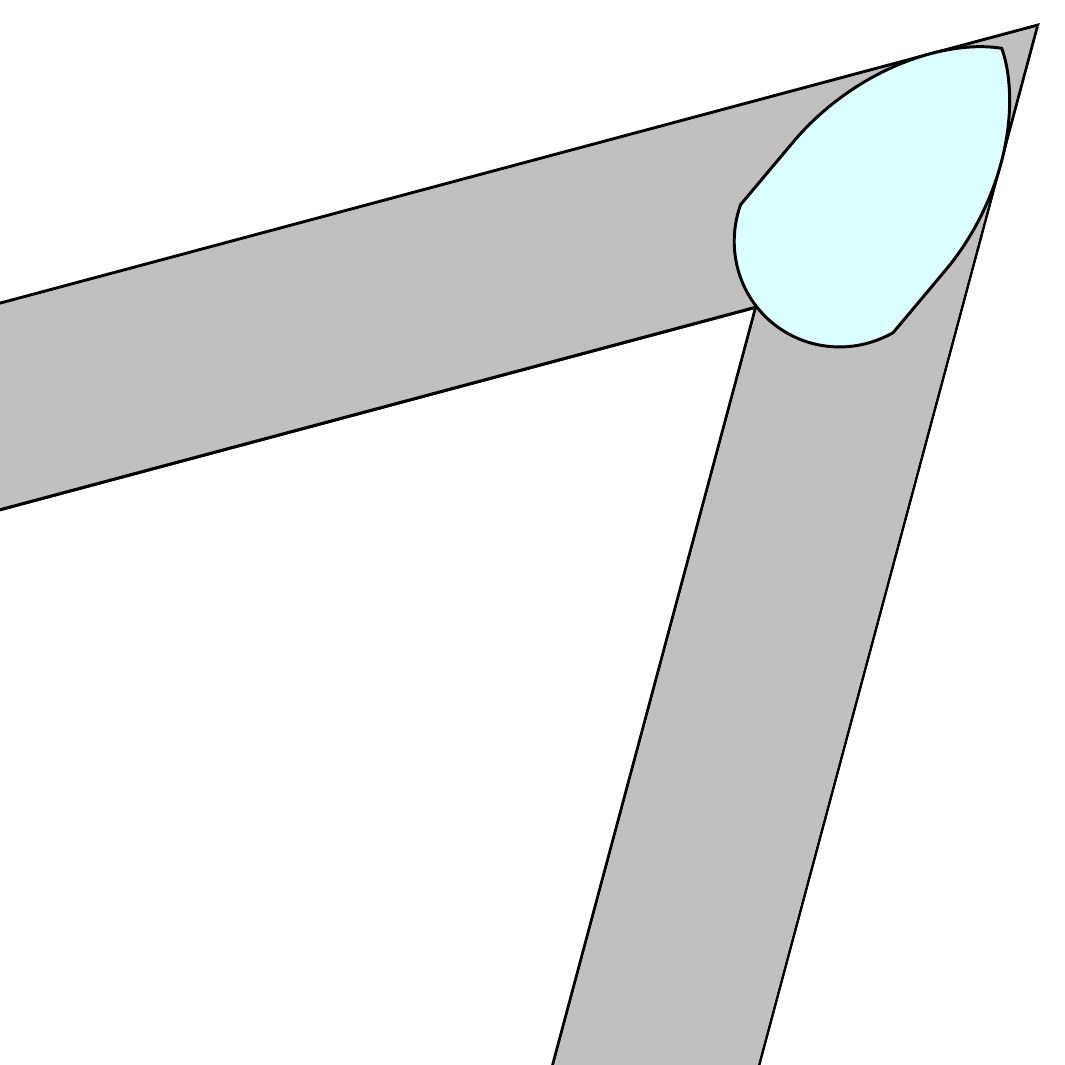}
\caption{$60^\circ$}
\end{subfigure}

\caption{Two motion patterns for three corridor angles}
\label{fig:angle_compare}

\end{figure}

We propose a numerical approach based on dynamical systems to study the moving sofa problem. Inspired by a physical model of gas-pressure driving, the method transforms the geometric optimization problem into a rigid-body dynamical system. This framework allows us to systematically explore locally optimal sofa shapes for corridors with arbitrary angles, extending the classical right-angle setting.

In contrast to existing numerical approaches \cite{Maruyama1973,Gibbs2014,Leng2024}, the direction of evolution at each step in our method is directly determined by the dynamics, requiring no costly search or repeated evaluations of the area. Numerical experiments reproduce the known Gerver's sofa result and demonstrate that the method effectively explores feasible shapes across a wide range of corridor angles.

Our numerical experiments reveal an interesting phenomenon: two distinct motion patterns compete depending on the corridor angle. As illustrated in Figure~\ref{fig:angle_compare}, the dominant pattern changes at a critical angle of approximately $43.327^\circ$. At this transition point the locally maximal area is about $1.8674$. This suggests a \textbf { phase-transition-like } behavior in the optimal motion strategy as the corridor angle varies.

\begin{conjecture}[Phase transition of optimal motion patterns]
There exists a critical corridor angle $\theta_c$ at which the dominant optimal motion pattern switches. Numerically we observe
\[
\theta_c \approx 43.327^\circ ,
\qquad
A_{\max}(\theta_c) \approx 1.8674 .
\]
\end{conjecture}

The numerical results supporting Conjecture~1 are illustrated in Figure~\ref{fig:angle_compare}. 
For smaller corridor angles (e.g., $30^\circ$), the configuration shown in the second row produces larger feasible shapes and therefore dominates the optimization process. 
For larger angles (e.g., $60^\circ$), the configuration shown in the first row becomes more efficient and yields a larger area. 
Near $\theta \approx 43.327^\circ$, the two patterns produce nearly identical maximal areas, indicating a transition between the two regimes.

Geometrically, the two configurations correspond to different global motion strategies of the rigid body. 
The motion in the first row represents a complete clockwise rotation accompanied by translation, while the motion in the second row represents a complete counterclockwise rotation with translation. 
Our numerical evidence suggests that before the transition ($\theta < \theta_c$) the counterclockwise rotation--translation pattern is more efficient, whereas after the transition ($\theta > \theta_c$) the clockwise rotation--translation pattern becomes dominant. 
The data presented in Section~3 further support the phase-transition-like behavior described in Conjecture~1.

\section{Method}

In his study of the moving sofa problem, Gerver proposed the balanced polygon condition, suggesting that the boundary of an optimal configuration must satisfy a kind of equivalent pressure balance condition \cite{Gerver1992}. This condition indicates that the optimal shape can be understood as a static equilibrium state under pressures acting in several directions.

The basic idea of this paper is to interpret this static equilibrium condition as a degenerate rigid-body dynamical process. Specifically, we impose a uniform normal pressure on the boundary of the feasible region and construct a corresponding equivalent potential energy, thereby transforming the geometric optimization problem into a gradient flow system in the configuration space.

\subsection{Algorithm}

Below we establish the correspondence between boundary pressure and the shape derivative of the area functional, and from this derive the corresponding degenerate gradient flow dynamics to characterize local extremal configurations of the geometric intersection area. For convenience, we adopt the viewpoint of fixing the sofa and moving the corridor, taking the intersection of the corridors at all times as the sofa area.

Let

\[
q \in \mathbb{R}^k
\]

be a finite-dimensional configuration variable describing the relative positions and orientations of several corridors. For each admissible configuration $q$, its geometry determines a planar region

\[
\Omega(q) \subset \mathbb{R}^2 ,
\]

and we denote its area by

\[
A(q) := |\Omega(q)| .
\]

The concern of this paper is to characterize configurations $q$ that yield a local maximum of $A(q)$. Throughout the derivation, we do not require $A(q)$ to have an explicit analytic expression, but we assume that $A(q)$ is piecewise differentiable with respect to $q$.

Consider a virtual displacement on the boundary $\partial \Omega$:

\[
\delta \mathbf{x} : \partial \Omega \to \mathbb{R}^2 .
\]

If a constant normal pressure $p>0$ is applied on the boundary, the virtual work done by this pressure on the virtual displacement is defined as

\[
\delta W
=
\int_{\partial \Omega}
p \, \mathbf{n} \cdot \delta \mathbf{x} \, \mathrm{d}s ,
\tag{1}
\]

where $\mathbf{n}$ denotes the outward unit normal vector. This expression depends only on geometric and kinematic facts: pressure does work only on the normal component of the boundary displacement, independent of any specific physical model.

Introduce a scalar functional $\Pi(\Omega)$ defined on planar regions, and require that its first-order shape variation satisfies

\[
\delta \Pi(\Omega) = - \delta W .
\tag{2}
\]

Assume that the shape potential energy has a volume integral representation

\[
\Pi(\Omega)
=
- \int_{\Omega} f(x) \, \mathrm{d}x ,
\]

where the function $f$ is independent of the shape of the region. According to the structure theorem for Hadamard shape derivatives \cite{hadamard1923}, its first-order shape variation is

\[
\delta \Pi
=
- \int_{\partial \Omega}
f(x)\, \mathbf{n}\cdot \delta \mathbf{x} \, \mathrm{d}s .
\tag{3}
\]

Comparing (3) with (1)--(2), we see that if for every admissible boundary virtual displacement we have

\[
\delta \Pi = - \delta W ,
\]

then necessarily

\[
f(x) \equiv p .
\]

Thus, the shape potential can be taken as

\[
\Pi(\Omega) = - p |\Omega| .
\]

This potential is not interpreted as the internal energy of a real physical system, but rather as an equivalent variational potential whose sole purpose is to generate the prescribed boundary pressure through variation.

Since the region $\Omega$ is determined by the configuration variable $q$, the shape potential can be written as

\[
\Pi(q) := \Pi(\Omega(q)) = - p A(q) .
\]

A virtual variation $\delta q$ of the configuration induces a displacement on the boundary:

\[
\delta \mathbf{x}
=
\frac{\partial \mathbf{x}}{\partial q} \, \delta q .
\]

Substituting into the shape variation formula yields

\[
\delta \Pi
=
- \int_{\partial \Omega(q)}
p \,
\mathbf{n}\cdot
\frac{\partial \mathbf{x}}{\partial q}
\, \delta q \, \mathrm{d}s .
\]

Therefore, we can define the corresponding generalized force in configuration space as

\[
F(q)
:=
- \nabla_q \Pi(q)
=
p \, \nabla_q A(q) .
\]

This definition does not rely on an explicit form of $A(q)$; it is entirely induced by the geometric variation of the boundary.

\begin{figure}[htbp]
\centering
\includegraphics[width=1\textwidth]{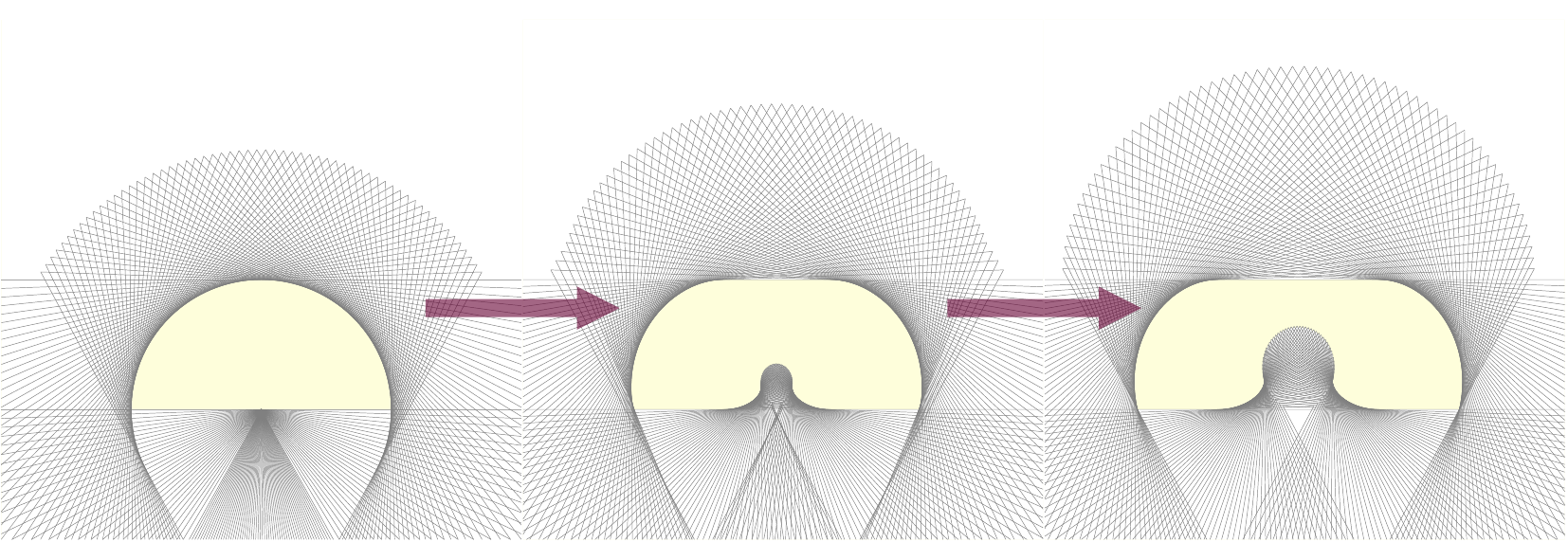}
\caption{Evolution process for a 60° corridor configuration}
\label{fig:tu}
\end{figure}

In the limit where inertia is neglected and only dissipative effects are retained, for configuration regions where the area function is smooth, the evolution of the configuration variables satisfies

\[
\gamma \dot q
=
- \nabla_q \Pi(q),
\qquad \gamma>0 ,
\]

where $\gamma$ is the damping coefficient. Substituting the specific expression for the shape potential gives

\[
\dot q
=
\frac{p}{\gamma} \nabla_q A(q) .
\tag{4}
\]

Along any smooth trajectory $q(t)$, we have

\[
\frac{\mathrm{d}}{\mathrm{d}t} \Pi(q(t))
=
- \frac{1}{\gamma}
\bigl\lvert \nabla_q \Pi(q(t)) \bigr\rvert^2
\le 0 ,
\]

equivalently,

\[
\frac{\mathrm{d}}{\mathrm{d}t} A(q(t))
=
\frac{p}{\gamma}
\bigl\lvert \nabla_q A(q(t)) \bigr\rvert^2
\ge 0 .
\]

Thus, the shape potential $\Pi$ is a Lyapunov function for this dynamical system, and the area of the region is non-decreasing along trajectories; Figure~\ref{fig:tu} illustrates the evolution process for a 60° corridor configuration.

At configuration points where the area function is smooth, an equilibrium configuration $q^\ast$ satisfies

\[
\nabla_q A(q^\ast) = 0 .
\]

If the Hessian matrix $\nabla_q^2 A(q^\ast)$ is negative semidefinite, then this equilibrium is stable in the dynamical sense and corresponds to a local maximum of the area functional. If the Hessian is indefinite, the equilibrium is an unstable saddle point. This paper does not address questions of existence or uniqueness of global optimal solutions.

The above derivation shows that a geometric optimization problem can be transformed into a degenerate gradient flow system via the potential structure induced by boundary variation, without explicitly constructing the gradient of the objective function. In this framework, local extrema of geometric quantities are naturally characterized by dynamical stability, rather than obtained by directly solving an optimization problem.

\subsection{Geometry}

\begin{figure}[htbp]
\centering
\includegraphics[width=0.8\textwidth]{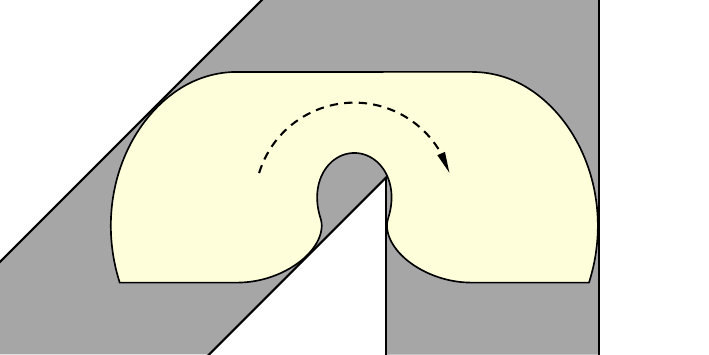}
\caption{Motion pattern}
\label{fig:ro1}
\end{figure}

We only study sofas that rotate monotonically and fully traverse the corridor; this appears to be a highly plausible necessary condition for an optimal sofa. This section discusses the motion pattern shown in Figure~\ref{fig:ro1}, namely that the sofa rotates clockwise during its motion and completes a total rotation equal to the supplement of the corridor angle. As will be demonstrated in Section~2.4, specifying the motion pattern is essential, because different motion patterns lead to completely different locally optimal shapes.

\begin{figure}[htbp]
\centering
\includegraphics[width=0.8\textwidth]{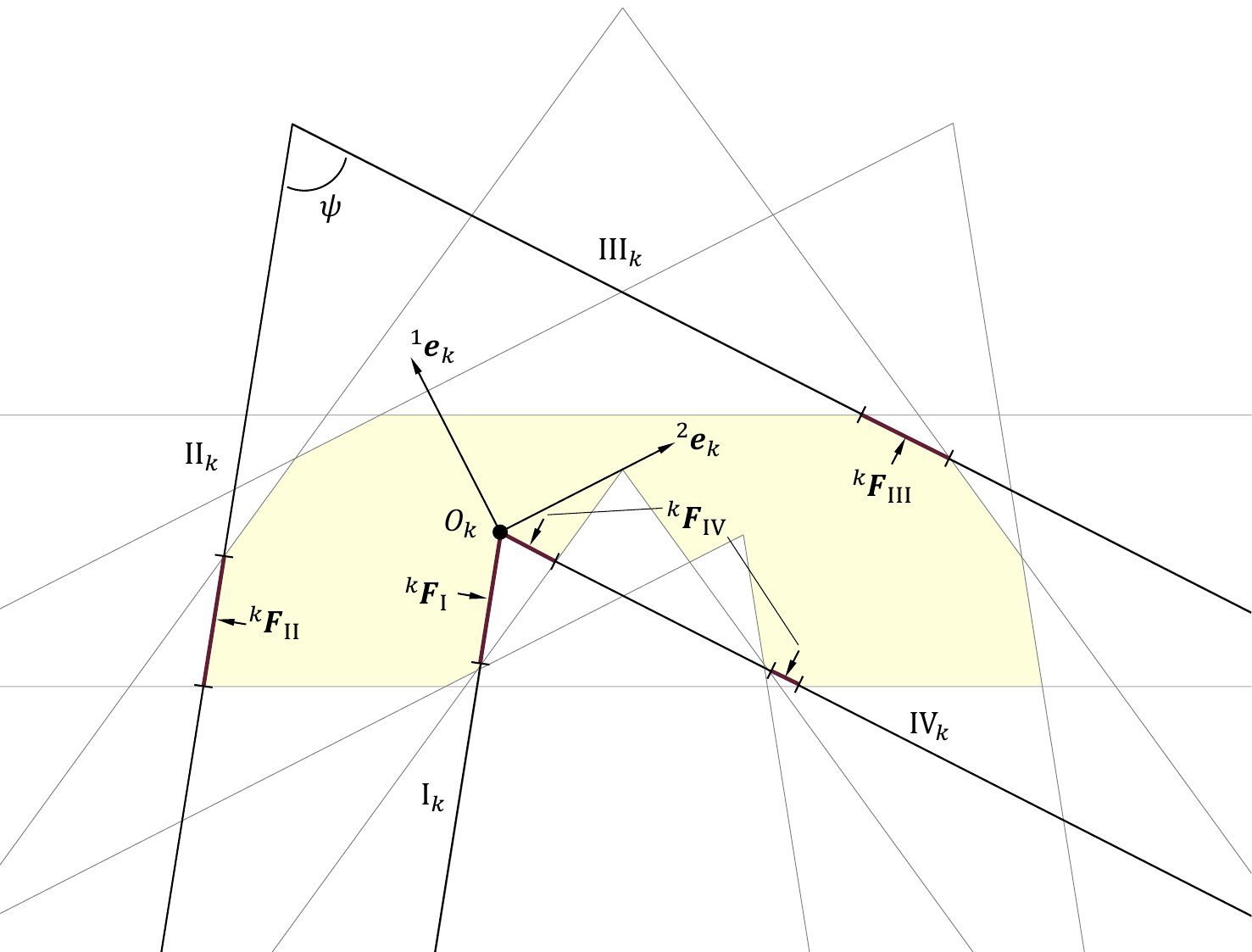}
\caption{Geometry of the corridor}
\label{fig:shape1}
\end{figure}

For numerical computation, we discretize the continuous rotation process into a finite number of corridor orientations. The orientation angles of the corridor during rotation are assumed to be equally spaced, as shown in Figure~\ref{fig:shape1}, and are kept fixed throughout the computation. 
The geometric configuration of the corridor is parameterized by the reference point positions of $n$ corridors. From Gerver's analysis of local optimality conditions, the first and last corridors each contribute only one branch (the horizontal branch) to the sofa intersection; therefore we merge them into a single fixed horizontal region.

Let $\psi$ be the angle, and let $k=1,2,\dots,n$ be integers. Denote the $k$-th reference point by
\[
O_k(\xi_k,\eta_k).
\]

Define two families of reference direction vectors:
\[
{}^{1}\mathbf e_k=
\bigl(
\cos\alpha_k,
\sin\alpha_k
\bigr),
\qquad
{}^{2}\mathbf e_k=
\bigl(
\cos\beta_k,
\sin\beta_k
\bigr),
\]
where
\[
\alpha_k=\frac{\psi}{2}+\frac{k(\pi-\psi)}{n+1},
\qquad
\beta_k=\frac{\psi-\pi}{2}+\frac{k(\pi-\psi)}{n+1}.
\]

Denote the four sides of the corridor by $\mathrm I_k,\mathrm{II}_k,\mathrm{III}_k,\mathrm{IV}_k$, whose line equations are uniformly written as
\[
y={}^{k}K_\Omega x+{}^{k}Y_\Omega,
\qquad
\Omega\in\{\mathrm I,\mathrm{II},\mathrm{III},\mathrm{IV}\}.
\]

The slopes of the four sides are respectively
\[
{}^{k}K_{\mathrm I}
={}^kK_{\mathrm{II}}
=\tan\frac{k(\pi-\psi)}{n+1},
\]
\[
{}^{k}K_{\mathrm{III}}
={}^kK_{\mathrm{IV}}
=\tan\!\left(\psi+\frac{k(\pi-\psi)}{n+1}\right).
\]

The intercepts are respectively
\[
{}^{k}Y_{\mathrm I}
=\eta_k-\xi_k\,{}^{k}K_{\mathrm I},
\qquad
{}^{k}Y_{\mathrm{IV}}
=\eta_k-\xi_k\,{}^{k}K_{\mathrm{IV}},
\]
\[
{}^{k}Y_{\mathrm{II}}
=\eta_k-\xi_k\,{}^{k}K_{\mathrm{II}}
+\frac{
\sin\alpha_k-{}^{k}K_{\mathrm{II}}\cos\alpha_k
}{
\sin(\psi/2)
},
\]
\[
{}^{k}Y_{\mathrm{III}}
=\eta_k-\xi_k\,{}^{k}K_{\mathrm{III}}
+\frac{
\sin\alpha_k-{}^{k}K_{\mathrm{III}}\cos\alpha_k
}{
\sin(\psi/2)
}.
\]

The corresponding $x$-direction truncation conditions are
\[
x\,\mathrm{sgn}({}^{k}K_{\mathrm I})
\le
\xi_k\,\mathrm{sgn}({}^{k}K_{\mathrm I}),
\]
\[
x\,\mathrm{sgn}({}^{k}K_{\mathrm{II}})
\le
\left(
\xi_k+\frac{\cos\alpha_k}{\sin(\psi/2)}
\right)
\mathrm{sgn}({}^{k}K_{\mathrm{II}}),
\]
\[
x\,\mathrm{sgn}({}^{k}K_{\mathrm{III}})
\le
\left(
\xi_k+\frac{\cos\alpha_k}{\sin(\psi/2)}
\right)
\mathrm{sgn}({}^{k}K_{\mathrm{III}}),
\]
\[
x\,\mathrm{sgn}({}^{k}K_{\mathrm{IV}})
\le
\xi_k\,\mathrm{sgn}({}^{k}K_{\mathrm{IV}}).
\]

\subsection{Boundary Pressure}

To discretize the shape variation framework established in the previous section onto a concrete geometric configuration, we decompose the action of the boundary pressure onto each contact side, as shown in Figure~\ref{fig:shape1}. Let

\[
{}^{k}F_{\mathrm I},
\ {}^{k}F_{\mathrm{II}},
\ {}^{k}F_{\mathrm{III}},
\ {}^{k}F_{\mathrm{IV}}
\]

be the components of the resultant force produced by the uniform normal pressure on the corresponding boundary segments. The magnitude of each component equals the product of the pressure and the effective contact length of that side, and its direction is determined by the outward normal of the corresponding side.

Consequently, the total force acting on the $k$-th frame can be expressed as the vector sum of the pressure resultants on its four sides:

\[
{}^{k}\mathbf F
=
{}^{k}\mathbf F_{\mathrm I}
+
{}^{k}\mathbf F_{\mathrm{II}}
+
{}^{k}\mathbf F_{\mathrm{III}}
+
{}^{k}\mathbf F_{\mathrm{IV}}.
\]

In the numerical implementation, we do not explicitly compute the area gradient; instead, we approximate the generalized forces in configuration space by means of the aforementioned boundary resultants. According to the derivation in the previous section,

\[
{}^{k}\mathbf F
=
\nabla_{q_k} A(q),
\]

so each component force ${}^{k}F_{\Omega}$ can be viewed as the variational contribution of the area to the local displacement of that side. By superimposing the pressure resultants on all contact sides, we obtain a discrete approximation of the gradient of the area functional in configuration space. From this we construct a discrete gradient flow for the area functional, which automatically drives the configuration toward a state of pressure balance and area extremum.

The resultant force on the $k$-th frame in global coordinates is

\[
\begin{pmatrix}
{}^{k}F_x \\ {}^{k}F_y
\end{pmatrix}
=
\begin{pmatrix}
\cos\beta_k & -\sin\beta_k \\
\sin\beta_k & \cos\beta_k
\end{pmatrix}
\begin{pmatrix}
\cos(\psi/2)\bigl({}^{k}F_{\mathrm I}-{}^{k}F_{\mathrm{II}}+{}^{k}F_{\mathrm{III}}-{}^{k}F_{\mathrm{IV}}\bigr) \\
\sin(\psi/2)\bigl(-{}^{k}F_{\mathrm I}+{}^{k}F_{\mathrm{II}}+{}^{k}F_{\mathrm{III}}-{}^{k}F_{\mathrm{IV}}\bigr)
\end{pmatrix}.
\]

The corresponding degenerate gradient flow dynamical system is

\[
\begin{cases}
\dot{\xi}_k = \dfrac{1}{\gamma} \, {}^{k}F_x \\[6pt]
\dot{\eta}_k = \dfrac{1}{\gamma} \, {}^{k}F_y
\end{cases}, \quad k=1,2,\dots,n.
\]

\subsection{Intersection}

To compute the generalized forces in configuration space, we need to determine the effective contact length of each corridor boundary within the intersection region. This length equals the segment length of the corresponding boundary line within the intersection region of all corridors.

\begin{figure}[htbp]
\centering
\includegraphics[width=0.8\textwidth]{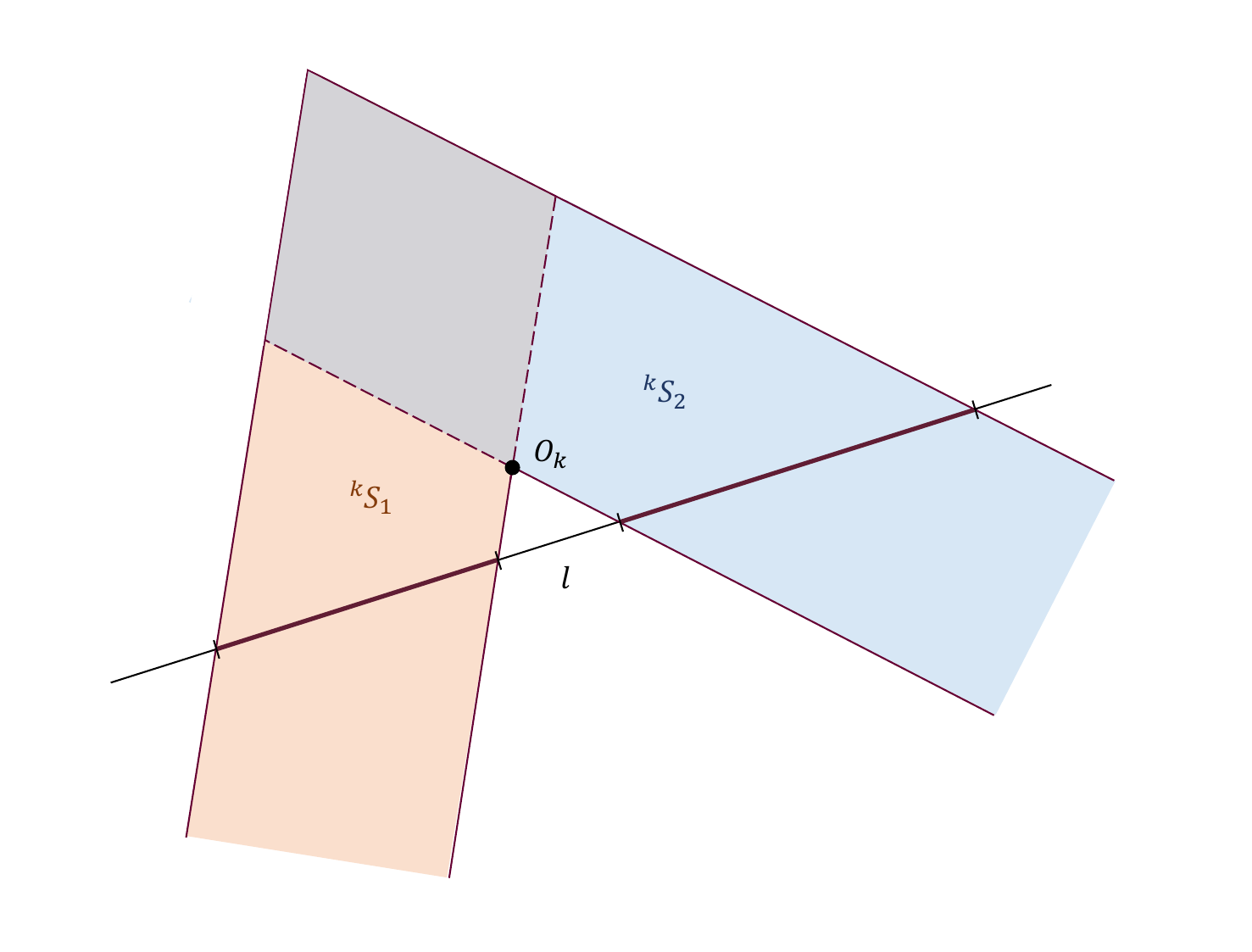}
\caption{Intersection length}
\label{fig:shape2}
\end{figure}

The region of a line $l$ inside the corridor intersection is

\[
\Bigg(
\bigcap_{k=1}^{n}
\Big(
\big({}^{k}S_1 \cap l\big)
\cup
\big({}^{k}S_2 \cap l\big)
\Big)
\Bigg)
\cap S
\tag{5}
\]

where ${}^{k}S_1$ and ${}^{k}S_2$ are the two branches of the $k$-th corridor, the two branches have an overlapping grey region (see Figure~\ref{fig:shape2}), and $S$ is the horizontal region.

Let $l$ be

\[
\Omega_m, \quad m = 1,2,\dots,n
\]

We compute separately the $x$-range of the line $\Omega_m$ under the truncation by each side of every corridor:

\[
x({}^{m}K_{\Omega} - {}^{k}K_{\text{I}}) \cdot \mathrm{sgn}({}^{k}K_{\text{I}}) 
\ge 
\left( {}^{k}Y_{\text{I}} - {}^{m}Y_{\Omega} \right) \cdot \mathrm{sgn}({}^{k}K_{\text{I}})
\tag{6}
\]

\[
x({}^{m}K_{\Omega} - {}^{k}K_{\text{II}}) \cdot \mathrm{sgn}({}^{k}K_{\text{II}}) 
\le 
\left( {}^{k}Y_{\text{II}} - {}^{m}Y_{\Omega} \right) \cdot \mathrm{sgn}({}^{k}K_{\text{II}})
\tag{7}
\]

\[
x({}^{m}K_{\Omega} - {}^{k}K_{\text{III}}) \cdot \mathrm{sgn}({}^{k}K_{\text{III}}) 
\ge 
\left( {}^{k}Y_{\text{III}} - {}^{m}Y_{\Omega} \right) \cdot \mathrm{sgn}({}^{k}K_{\text{III}})
\tag{8}
\]

\[
x({}^{m}K_{\Omega} - {}^{k}K_{\text{IV}}) \cdot \mathrm{sgn}({}^{k}K_{\text{IV}}) 
\le 
\left( {}^{k}Y_{\text{IV}} - {}^{m}Y_{\Omega} \right) \cdot \mathrm{sgn}({}^{k}K_{\text{IV}})
\tag{9}
\]

The $x$-range of the line $\Omega_m$ under the truncation by the horizontal region is

\[
{}^{m}K_{\Omega} \, x \ge - {}^{m}Y_{\Omega}, 
\quad
{}^{m}K_{\Omega} \, x \le 1 - {}^{m}Y_{\Omega}
\]

The $x$-range defined by the line $\Omega_m$ itself is

\[
\Phi_{\Omega} =
\begin{cases}
\{x: x \cdot \mathrm{sgn}({}^{k}K_{\mathrm I}) \le \xi_k \cdot \mathrm{sgn}({}^{k}K_{\mathrm I}) \}, & \Omega = \text{I} \\[0.5em]
\Big\{ x: x \cdot \mathrm{sgn}({}^{k}K_{\mathrm{II}}) \le \big( \xi_k + \frac{\cos\alpha_k}{\sin(\psi/2)} \big) \cdot \mathrm{sgn}({}^{k}K_{\mathrm{II}}) \Big\}, & \Omega = \text{II} \\[0.5em]
\Big\{ x: x \cdot \mathrm{sgn}({}^{k}K_{\mathrm{III}}) \le \big( \xi_k + \frac{\cos\alpha_k}{\sin(\psi/2)} \big) \cdot \mathrm{sgn}({}^{k}K_{\mathrm{III}}) \Big\}, & \Omega = \text{III} \\[0.5em]
\{x: x \cdot \mathrm{sgn}({}^{k}K_{\mathrm{IV}}) \le \xi_k \cdot \mathrm{sgn}({}^{k}K_{\mathrm{IV}}) \}, & \Omega = \text{IV}
\end{cases}
\]

Hence

\[
\{x: {}^{k}S_1 \cap \Omega_m\} = (6) \cap (7) \cap (8) \cap \Phi_{\Omega}
\]
\[
\{x: {}^{k}S_2 \cap \Omega_m\} = (7) \cap (8) \cap (9) \cap \Phi_{\Omega}
\]

From equation (5) we obtain the $x$-interval length of the line $\Omega_m$ within the intersection region, denoted by ${}^{m}L_{\Omega}$. Given the slope of the line, the corresponding actual segment length is

\[
{}^{m}b_{\Omega} = {}^{m}L_{\Omega} \sqrt{1 + \left({}^{m}K_{\Omega}\right)^2}
\]

The force is then computed as

\[
{}^{k}F_{\Omega} = p \cdot {}^{k}b_{\Omega}
\]

\subsection{Alternative Motion Pattern}

\begin{figure}[htbp]
\centering
\includegraphics[width=0.8\textwidth]{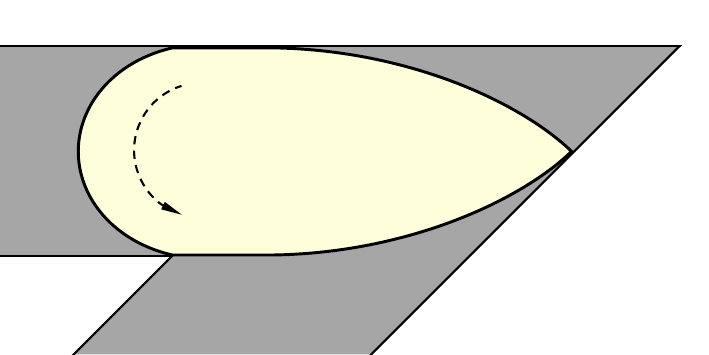}
\caption{Alternative motion pattern}
\label{fig:ro2}
\end{figure}

Different from the motion pattern described in Section 2.2, we may also consider an alternative motion pattern as shown in Figure~\ref{fig:ro2}, in which the sofa rotates counterclockwise overall while traversing the corridor, with a total rotation angle equal to the corridor angle. When the corridor angle is smaller than a certain value, this motion pattern yields a larger area compared to that of Section 2.2.

\begin{figure}[htbp]
\centering
\includegraphics[width=0.8\textwidth]{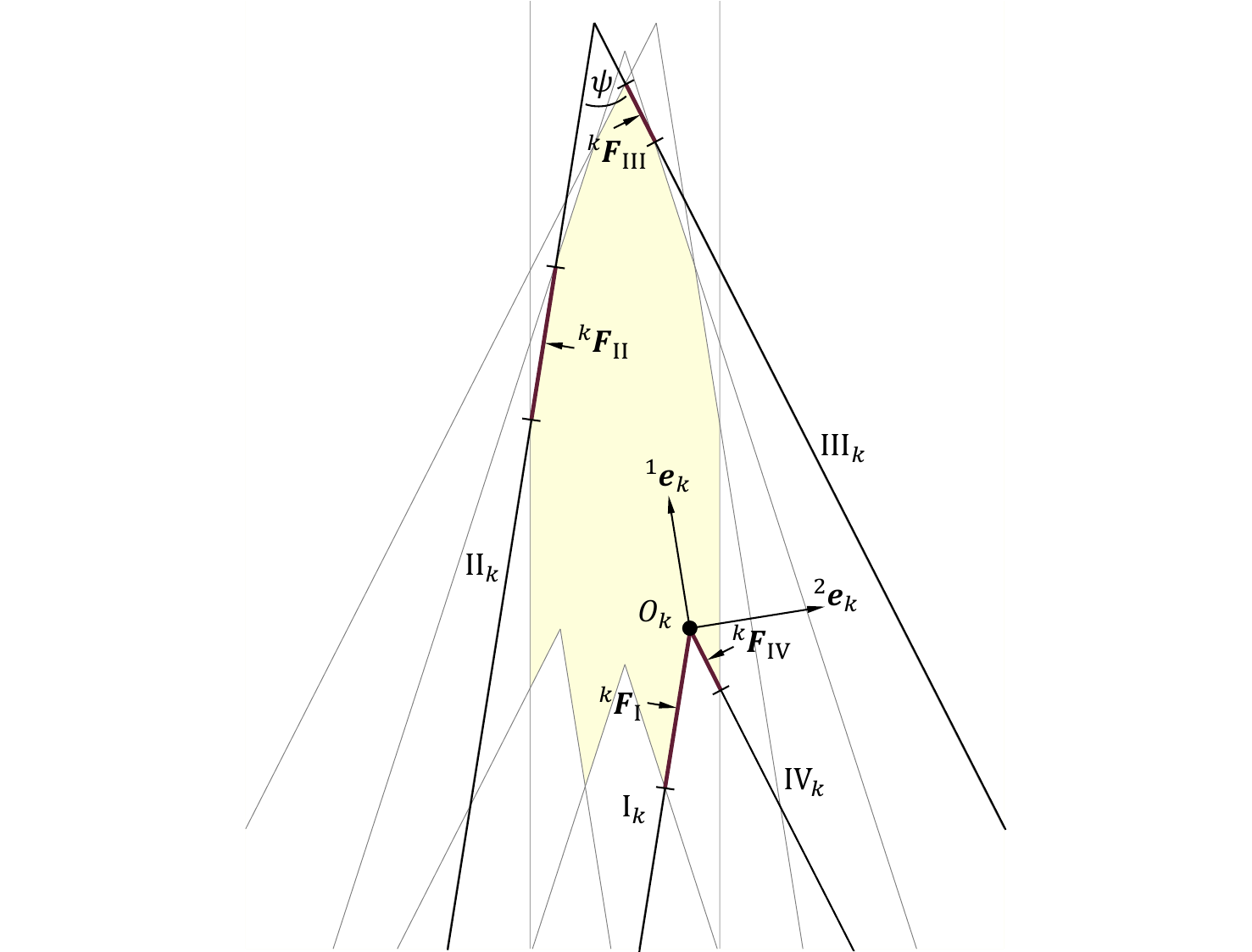}
\caption{Geometry of the corridor}
\label{fig:shape3}
\end{figure}

Using a discretization approach similar to that of Section 2.2, we discretize the continuous rotation process into $n$ equally spaced rotation angles, thereby obtaining $n$ corridors with evenly distributed angles and a fixed vertical region, as shown in Figure~\ref{fig:shape3}.

Let $\psi$ be the angle, and let $k=1,2,\dots,n$ be integers. Denote the $k$-th reference point by
\[
O_k(\xi_k,\eta_k).
\]

Define two families of reference direction vectors:
\[
{}^{1}\mathbf e_k=
\bigl(
\cos\alpha_k,
\sin\alpha_k
\bigr),
\qquad
{}^{2}\mathbf e_k=
\bigl(
\cos\beta_k,
\sin\beta_k
\bigr),
\]

where
\[
\alpha_k=\frac{\pi-\psi}{2}+\frac{k\psi}{n+1},
\qquad
\beta_k=-\frac{\psi}{2}+\frac{k\psi}{n+1}.
\]

Again denote the four sides by $\mathrm I_k,\mathrm{II}_k,\mathrm{III}_k,\mathrm{IV}_k$, whose line equations are uniformly written as
\[
y={}^{k}K_\Omega x+{}^{k}Y_\Omega,
\qquad
\Omega\in\{\mathrm I,\mathrm{II},\mathrm{III},\mathrm{IV}\}.
\]

The slopes of the four sides are
\[
{}^{k}K_{\mathrm I}
={}^kK_{\mathrm{II}}
=\tan\!\left(
\frac{\pi}{2}-\psi+\frac{k\psi}{n+1}
\right),
\]
\[
{}^{k}K_{\mathrm{III}}
={}^kK_{\mathrm{IV}}
=\tan\!\left(
\frac{\pi}{2}+\frac{k\psi}{n+1}
\right).
\]

The intercepts are respectively
\[
{}^{k}Y_{\mathrm I}
=\eta_k-\xi_k\,{}^{k}K_{\mathrm I},
\qquad
{}^{k}Y_{\mathrm{IV}}
=\eta_k-\xi_k\,{}^{k}K_{\mathrm{IV}},
\]
\[
{}^{k}Y_{\mathrm{II}}
=\eta_k-\xi_k\,{}^{k}K_{\mathrm{II}}
+\frac{
\sin\alpha_k-{}^{k}K_{\mathrm{II}}\cos\alpha_k
}{
\sin(\psi/2)
},
\]
\[
{}^{k}Y_{\mathrm{III}}
=\eta_k-\xi_k\,{}^{k}K_{\mathrm{III}}
+\frac{
\sin\alpha_k-{}^{k}K_{\mathrm{III}}\cos\alpha_k
}{
\sin(\psi/2)
}.
\]

The corresponding $x$-direction truncation conditions are
\[
x\,\mathrm{sgn}({}^{k}K_{\mathrm I})
\le
\xi_k\,\mathrm{sgn}({}^{k}K_{\mathrm I}),
\]
\[
x\,\mathrm{sgn}({}^{k}K_{\mathrm{II}})
\le
\left(
\xi_k+\frac{\cos\alpha_k}{\sin(\psi/2)}
\right)
\mathrm{sgn}({}^{k}K_{\mathrm{II}}),
\]
\[
x\,\mathrm{sgn}({}^{k}K_{\mathrm{III}})
\le
\left(
\xi_k+\frac{\cos\alpha_k}{\sin(\psi/2)}
\right)
\mathrm{sgn}({}^{k}K_{\mathrm{III}}),
\]
\[
x\,\mathrm{sgn}({}^{k}K_{\mathrm{IV}})
\le
\xi_k\,\mathrm{sgn}({}^{k}K_{\mathrm{IV}}).
\]

These expressions are similar in form to those of the previous framework; the steps for interval intersection and force calculation can be directly adopted from the preceding derivation.

\section{Results}

We implemented the above algorithm in \texttt{C++} and conducted numerical simulations for various corridor angles and numbers of frames. The dynamical system was integrated using a variable step-size method. The sofa area was computed using a polygon intersection algorithm, and the jagged boundaries resulting from contacts with the interior corners of the corridor were replaced by line segments connecting adjacent interior corner points to obtain a smooth approximate boundary. The results for the first motion pattern are shown in Tables~1 and~2, and those for the second motion pattern in Table~3. In particular, for the $90^\circ$ corridor case, the results in Table~1 are in close agreement with the Gerver's sofa area given in \cite{Gerver1992}, validating the effectiveness of the algorithm.

\begin{table}[htbp]
\centering
\small
\caption{Numerical results for different numbers of corridors $N$ (I)}
\begin{tabular}{c ccccc}
\toprule
$N$ 
& $\frac{\pi}{6}$ 
& $\frac{\pi}{4}$ 
& $\frac{\pi}{3}$ 
& $\frac{5\pi}{12}$ 
& $\frac{\pi}{2}$ \\
\midrule
100 
& 1.8202478345 
& 1.8744654111 
& 1.9508140523 
& 2.0595207893 
& 2.2195816868 \\

200 
& 1.8202209452 
& 1.8744281499 
& 1.9507825169 
& 2.0594837477 
& 2.2195474521 \\

300 
& 1.8202198265 
& 1.8744249978 
& 1.9507694919 
& 2.0594694911 
& 2.2195395225 \\

400 
& 1.8202196398 
& 1.8744240774 
& 1.9507673547 
& 2.0594668835 
& 2.2195342795 \\

500 
& 1.8202196525 
& 1.8744234554 
& 1.9507667152 
& 2.0594661375 
& 2.2195316460 \\
\bottomrule
\end{tabular}
\end{table}

\begin{table}[htbp]
\centering
\small
\caption{Numerical results for different numbers of corridors $N$ (II)}
\begin{tabular}{c ccccc}
\toprule
$N$ 
& $\frac{7\pi}{12}$ 
& $\frac{2\pi}{3}$ 
& $\frac{3\pi}{4}$ 
& $\frac{5\pi}{6}$ \\
\midrule
100 
& 2.4702997170 
& 2.8965448732 
& 3.6789651310 
& 5.3336855716 \\

200 
& 2.4702183816 
& 2.8963726641 
& 3.6786948617 
& 5.3332430670 \\

300 
& 2.4701984762 
& 2.8963305657 
& 3.6786317125 
& 5.3331421908 \\

400 
& 2.4701955339 
& 2.8963199291 
& 3.6786148651 
& 5.3331144254 \\

500 
& 2.4701928737 
& 2.8963150621 
& 3.6786071163 
& 5.3331016404 \\
\bottomrule
\end{tabular}
\end{table}

\begin{table}[htbp]
\centering
\small
\caption{Numerical results for the alternative rotation pattern for different numbers of corridors $N$}
\begin{tabular}{c cccccc}
\toprule
$N$ 
& $\frac{\pi}{12}$ 
& $\frac{\pi}{6}$ 
& $\frac{\pi}{4}$ 
& $\frac{\pi}{3}$ 
& $\frac{5\pi}{12}$ 
& $\frac{\pi}{2}$ \\
\midrule
100 
& 5.20637716 
& 2.64098072 
& 1.80373392 
& 1.39995665 
& 1.17172781 
& 1.03538276 \\

200 
& 5.20644124 
& 2.64101351 
& 1.80375639 
& 1.39997377 
& 1.17174133 
& 1.03539467 \\

300 
& 5.20645224 
& 2.64101902 
& 1.80376005 
& 1.39997644 
& 1.17174323 
& 1.03539581 \\

400 
& 5.20645653 
& 2.64102119 
& 1.80376148 
& 1.39997745 
& 1.17174395 
& 1.03539633 \\

500 
& 5.20645863 
& 2.64102225 
& 1.80376219 
& 1.39997795 
& 1.17174431 
& 1.03539661 \\

600 
& 5.20645977 
& 2.64102283 
& 1.80376258 
& 1.39997823 
& 1.17174451 
& 1.03539677 \\
\bottomrule
\end{tabular}
\end{table}

\begin{figure}[H]
\centering
\includegraphics[width=1\textwidth]{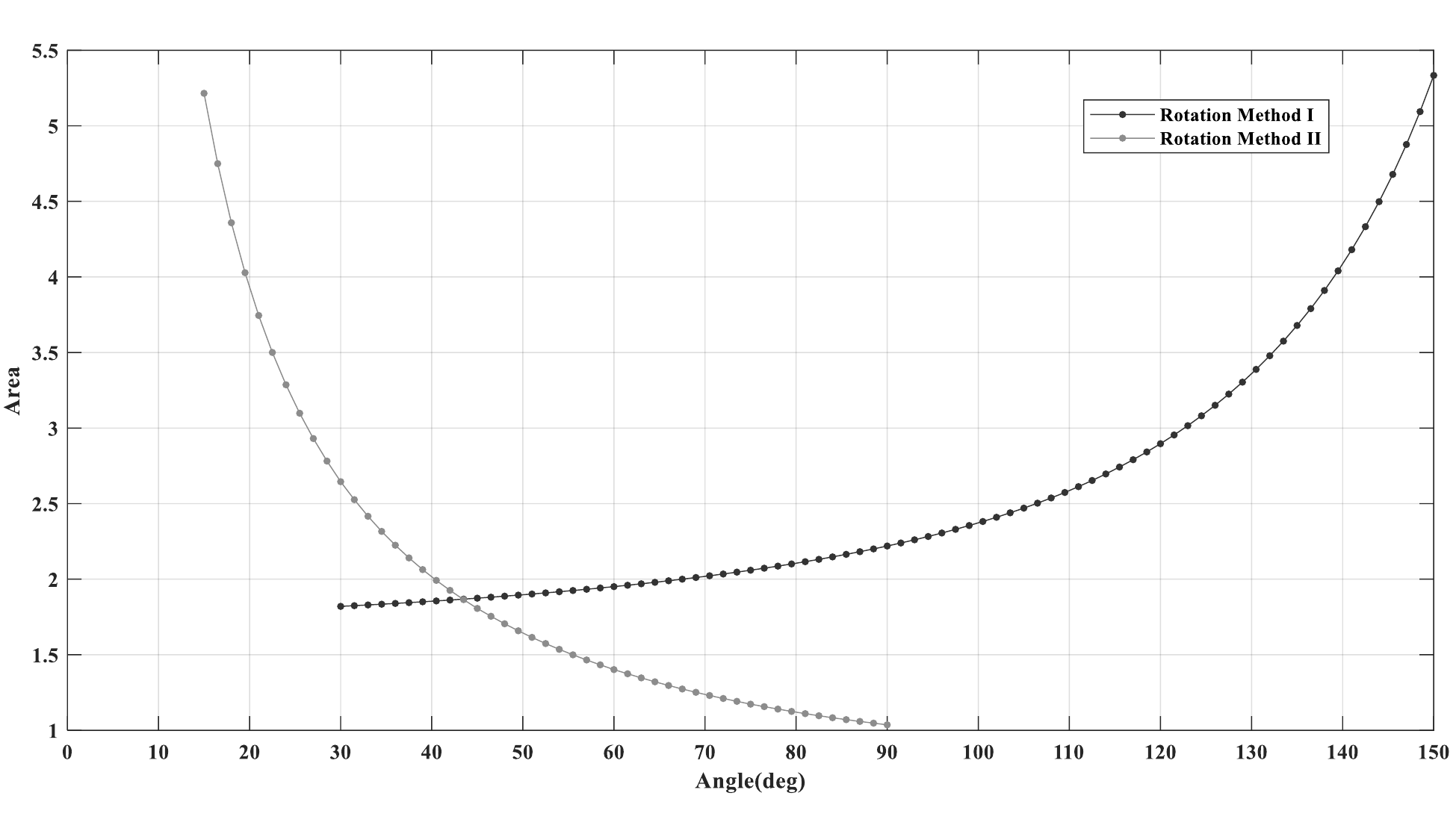}
\caption{Continuous variation of area with angle under the two patterns}
\label{fig:tu2}
\end{figure}

Figure~\ref{fig:tu2} shows the possible maximal sofa area for corridor angles from $15^\circ$ to $150^\circ$. For the first motion pattern, 100 corridors were used to compute locally optimal areas at 81 equally spaced angles from $30^\circ$ to $150^\circ$; for the second motion pattern, 300 corridors were used to compute locally optimal areas at 51 equally spaced angles from $15^\circ$ to $90^\circ$.

It can be observed that the two curves intersect between $\psi\in(43.327^\circ,\,43.328^\circ)$. For larger angles, the first motion pattern yields a larger area; for smaller angles, the second pattern performs better.

\begin{figure}[H]
\centering
\includegraphics[width=1\textwidth]{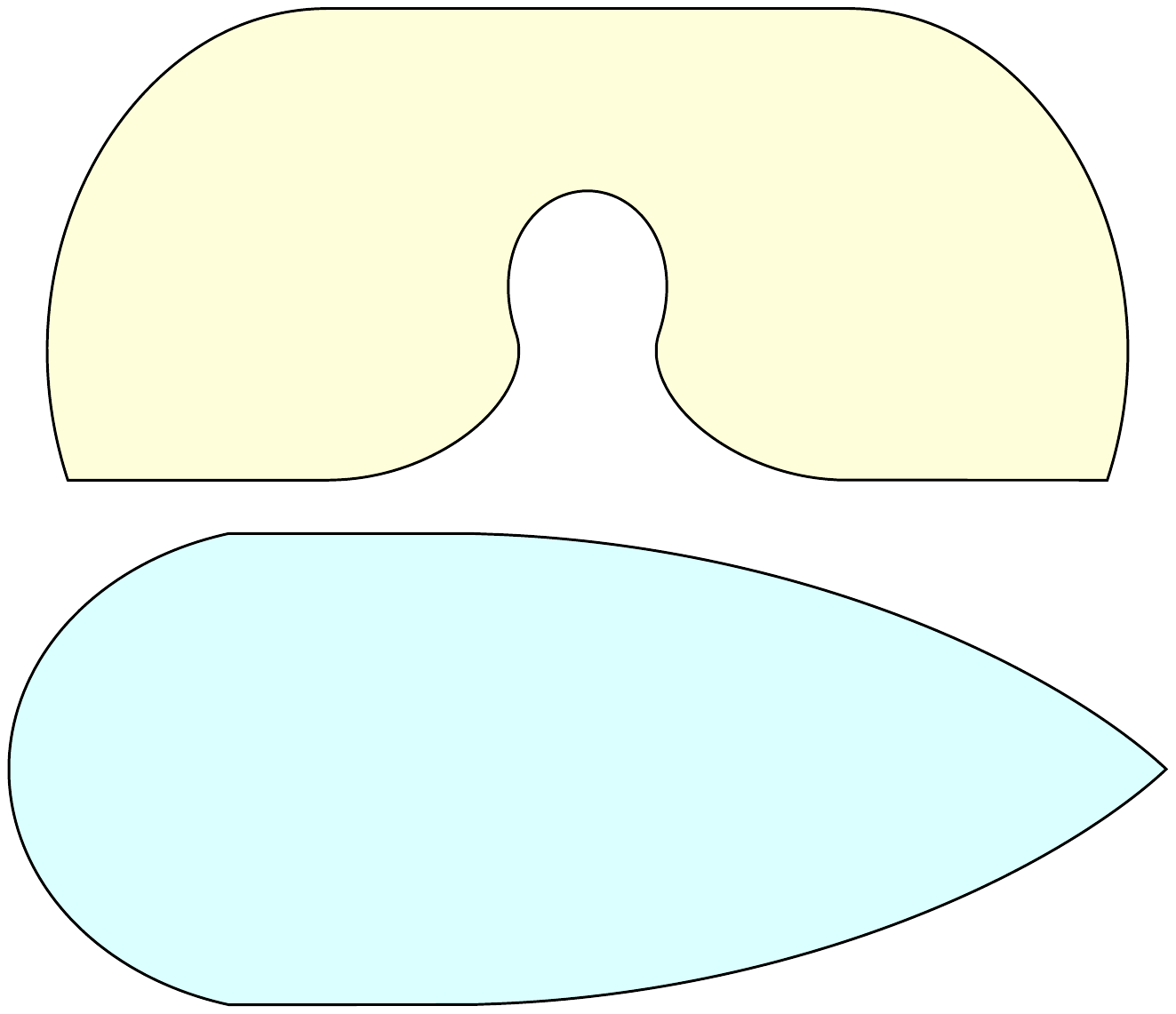}
\caption{Sofa shape for $\psi=43.327\ldots^\circ$, $Area=1.8674\ldots$}
\label{fig:43.327}
\end{figure}

To more precisely determine the intersection point of the two patterns, we performed higher-resolution numerical calculations at $\psi=43.327^\circ$ and $\psi=43.328^\circ$, with the results shown in Table~4. It can be seen that at $\psi=43.327^\circ$, the area obtained by the second motion pattern is slightly larger than that of the first; whereas at $\psi=43.328^\circ$, the first pattern becomes superior. This indicates that the area curves of the two patterns indeed intersect between these two angles. Figure~\ref{fig:43.327} shows the locally optimal sofa shape corresponding to $\psi=43.327\ldots^\circ$, with an area of approximately $1.8674\ldots$. This shape corresponds precisely to the critical region where the two motion patterns are nearly interchangeable.

\begin{table}[htbp]
  \centering
  \caption{Comparison of areas at $43.327^\circ$ and $43.328^\circ$}
  \small
  \begin{tabular}{c ccc}
    \toprule
    Angle & $N=100$ & $N=200$ & $N=300$ \\
    \midrule
    \multicolumn{4}{c}{\textbf{Pattern 1}} \\
    43.327° & 1.8674538445 & 1.8674175787 & 1.8674147097 \\
    43.328° & 1.8674579551 & 1.8674216923 & 1.8674188223 \\
    \midrule
    \multicolumn{4}{c}{\textbf{Pattern 2}} \\
    43.327° & 1.8674201401 & 1.8674418245 & 1.8674466806 \\
    43.328° & 1.8673805287 & 1.8674022128 & 1.8674070688 \\
    \bottomrule
  \end{tabular}
\end{table}

\begin{figure}[H]
\centering

\includegraphics[width=0.55\textwidth]{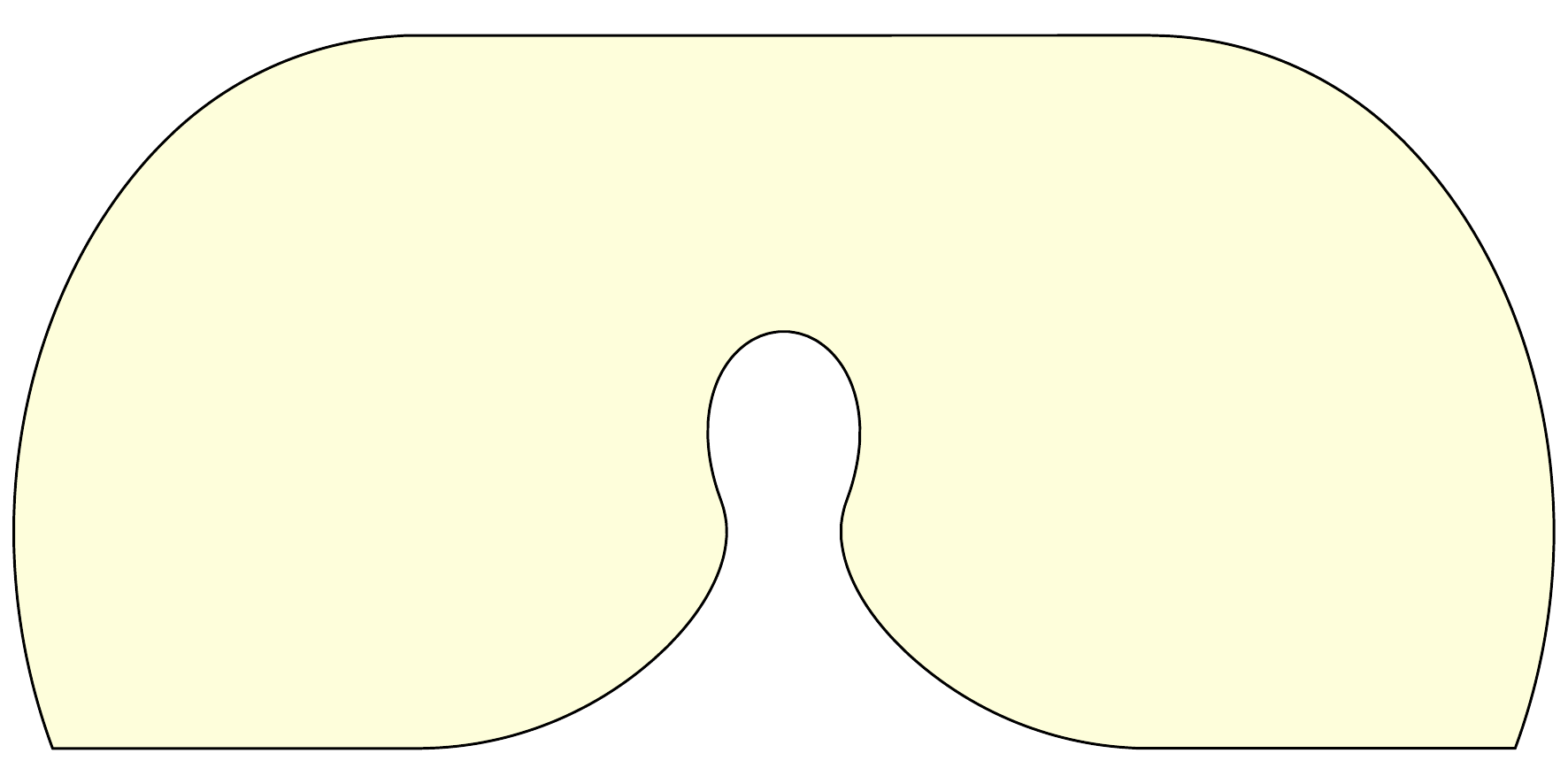}
\caption*{(a) $\psi=30^\circ$, $Area=1.8202...$}

\vspace{1.2em}

\includegraphics[width=0.55\textwidth]{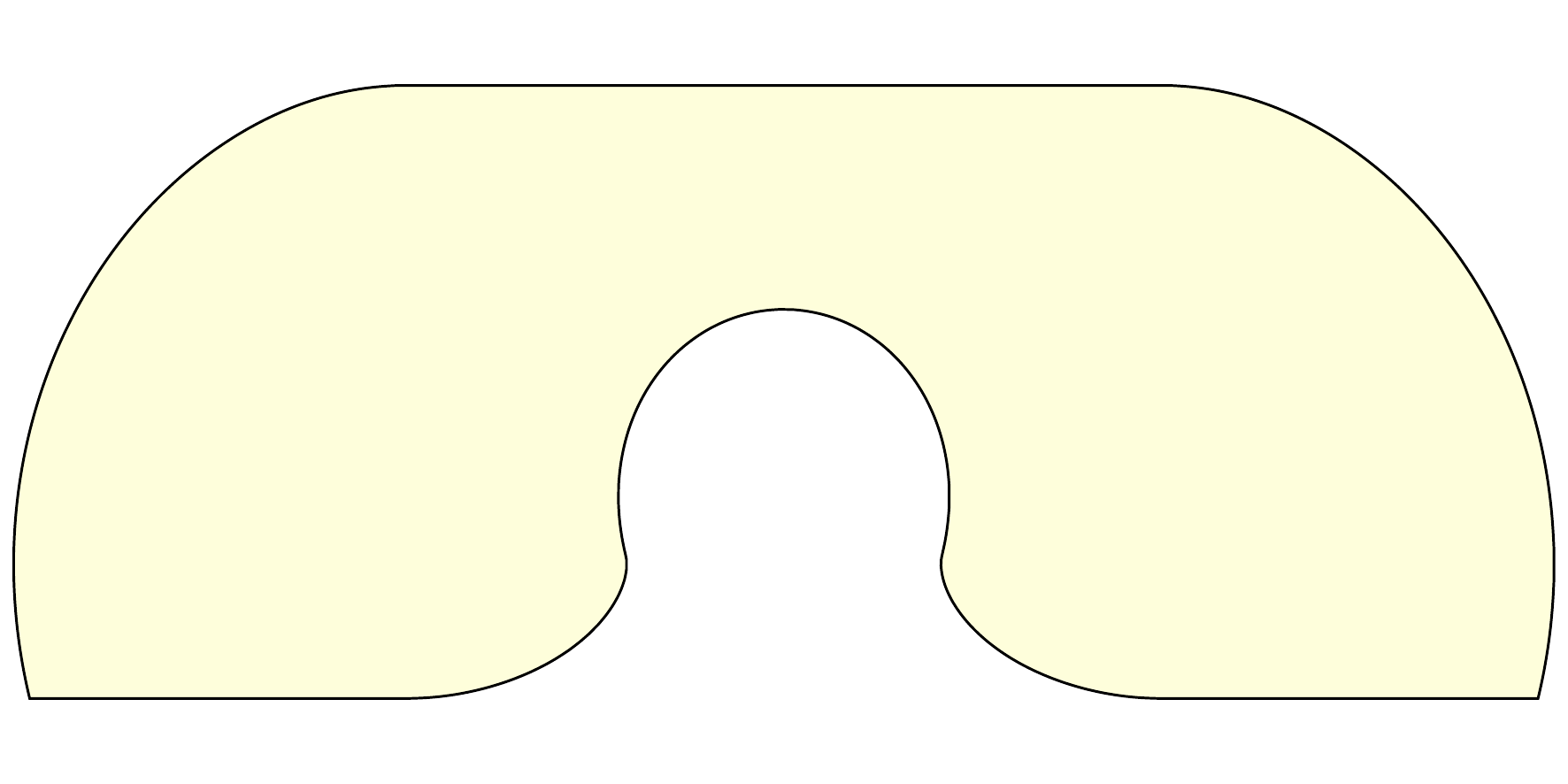}
\caption*{(b) $\psi=60^\circ$, $Area=1.9507...$}

\vspace{1.2em}

\includegraphics[width=0.55\textwidth]{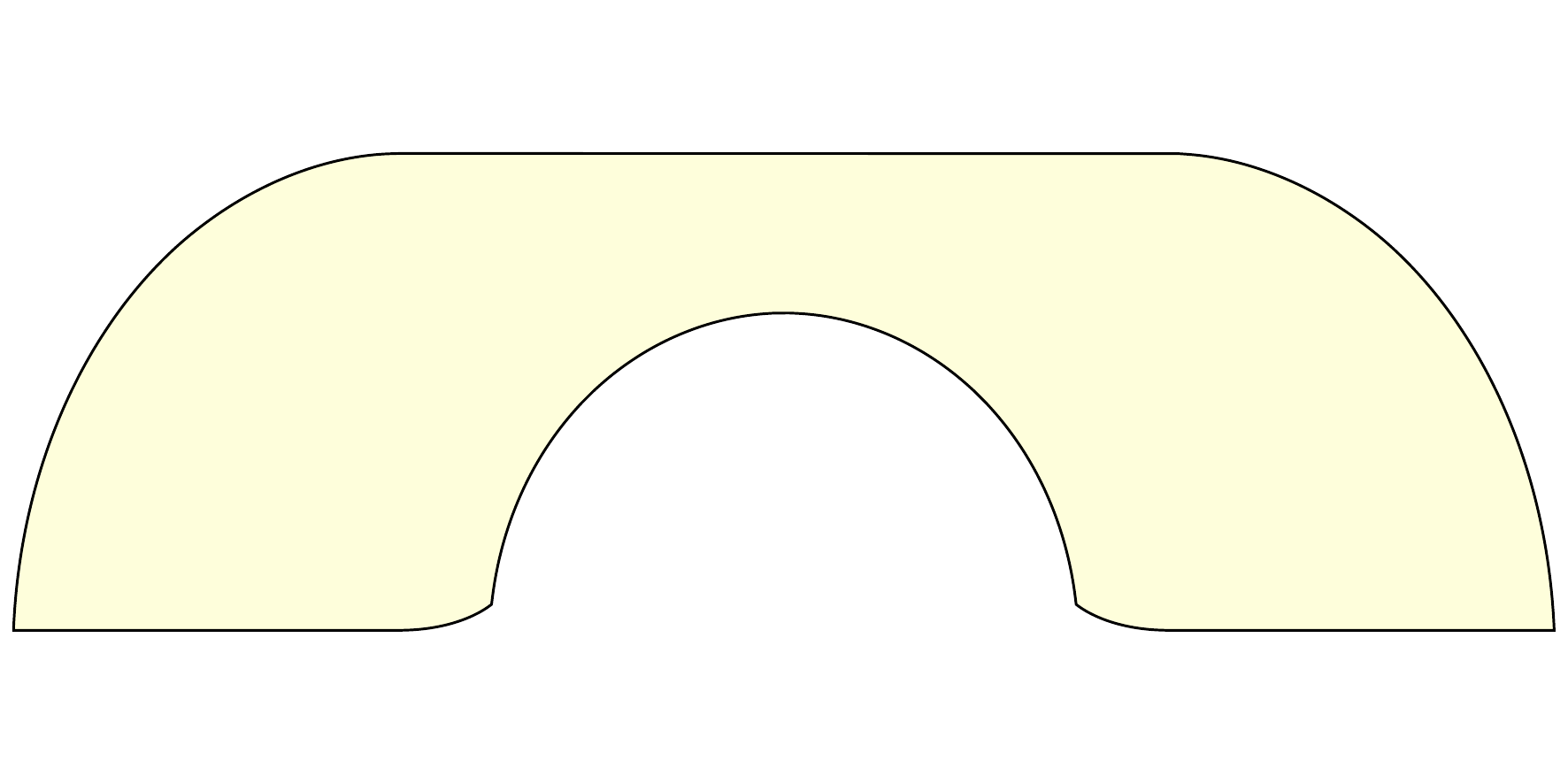}
\caption*{(c) $\psi=90^\circ$, $Area=2.2195...$}

\vspace{1.2em}

\includegraphics[width=0.55\textwidth]{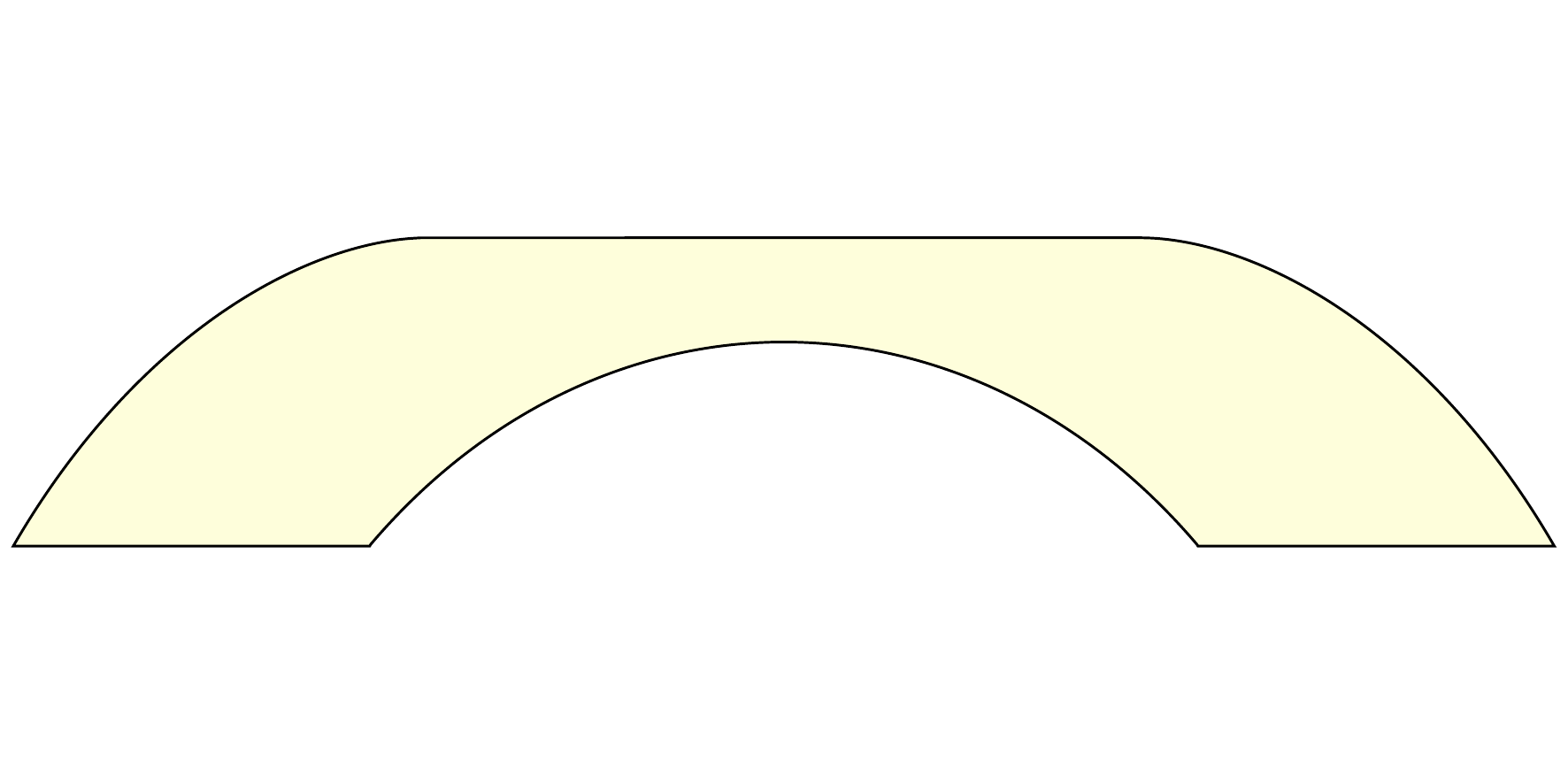}
\caption*{(d) $\psi=120^\circ$, $Area=2.8963...$}

\caption{Representative shapes for Pattern~1}

\label{fig:first}

\end{figure}

\begin{figure}[H]
\centering

\begin{subfigure}{0.34\textwidth}
\centering
\includegraphics[width=\linewidth]{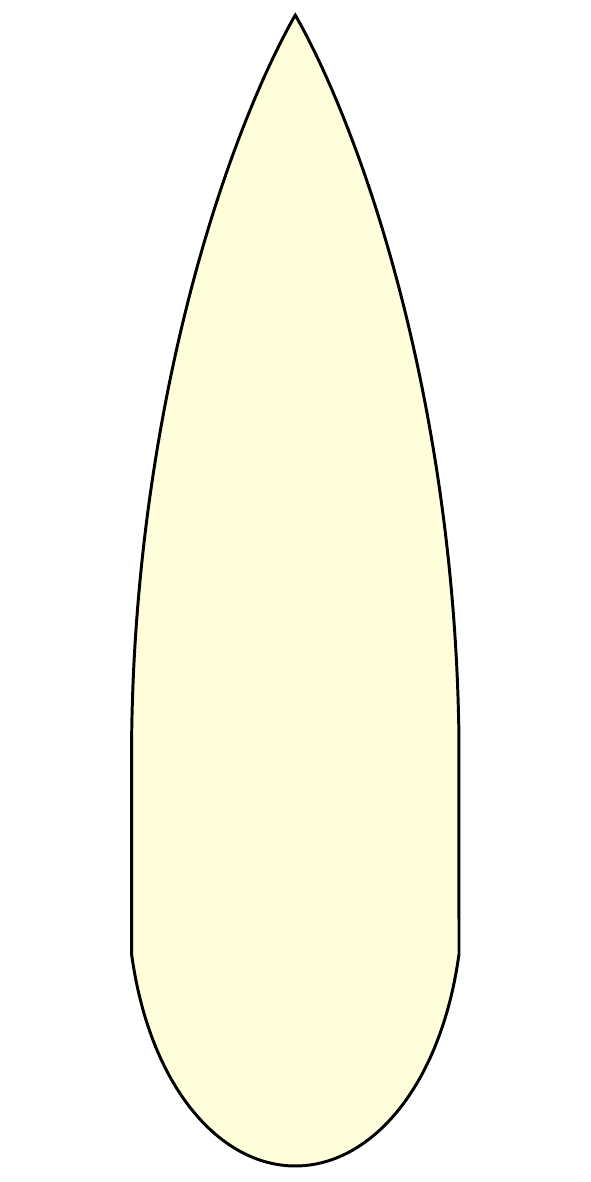}
\caption{$\psi=30^\circ$, $Area=2.6410...$}
\end{subfigure}
\hfill
\begin{subfigure}{0.34\textwidth}
\centering
\includegraphics[width=\linewidth]{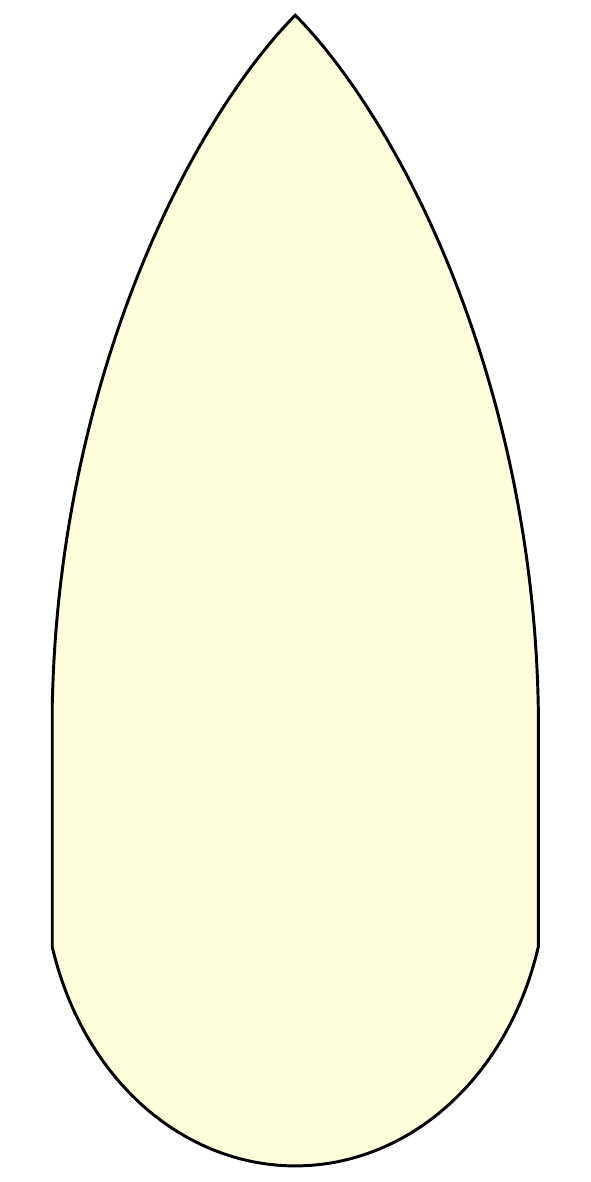}
\caption{$\psi=45^\circ$, $Area=1.8037...$}
\end{subfigure}

\vspace{1em}

\begin{subfigure}{0.34\textwidth}
\centering
\includegraphics[width=\linewidth]{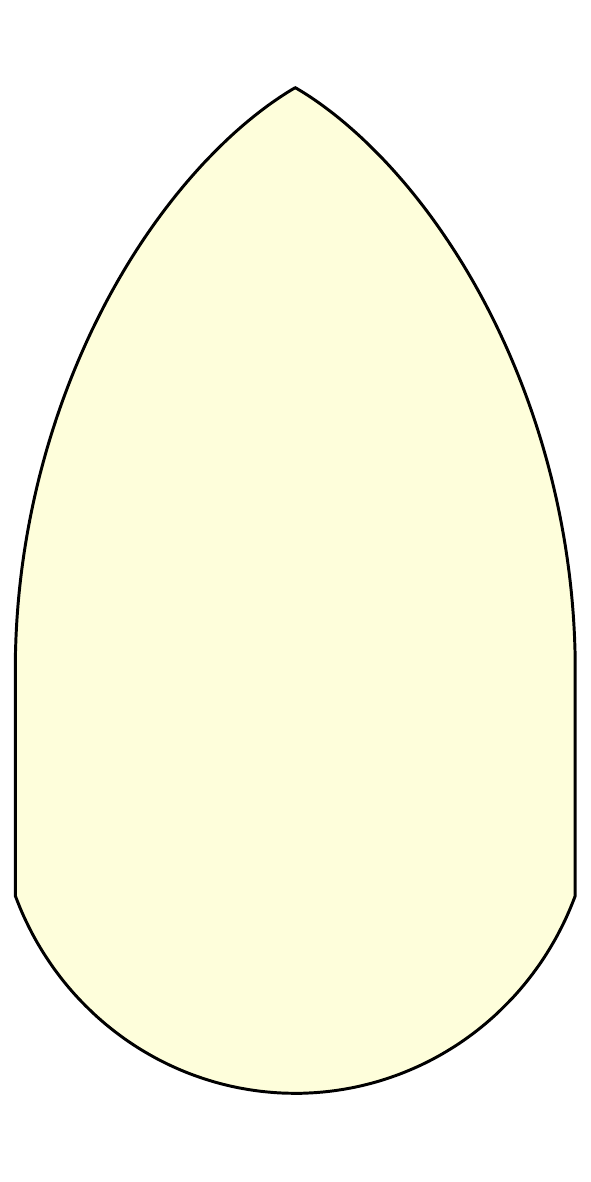}
\caption{$\psi=60^\circ$, $Area=1.3999...$}
\end{subfigure}
\hfill
\begin{subfigure}{0.34\textwidth}
\centering
\includegraphics[width=\linewidth]{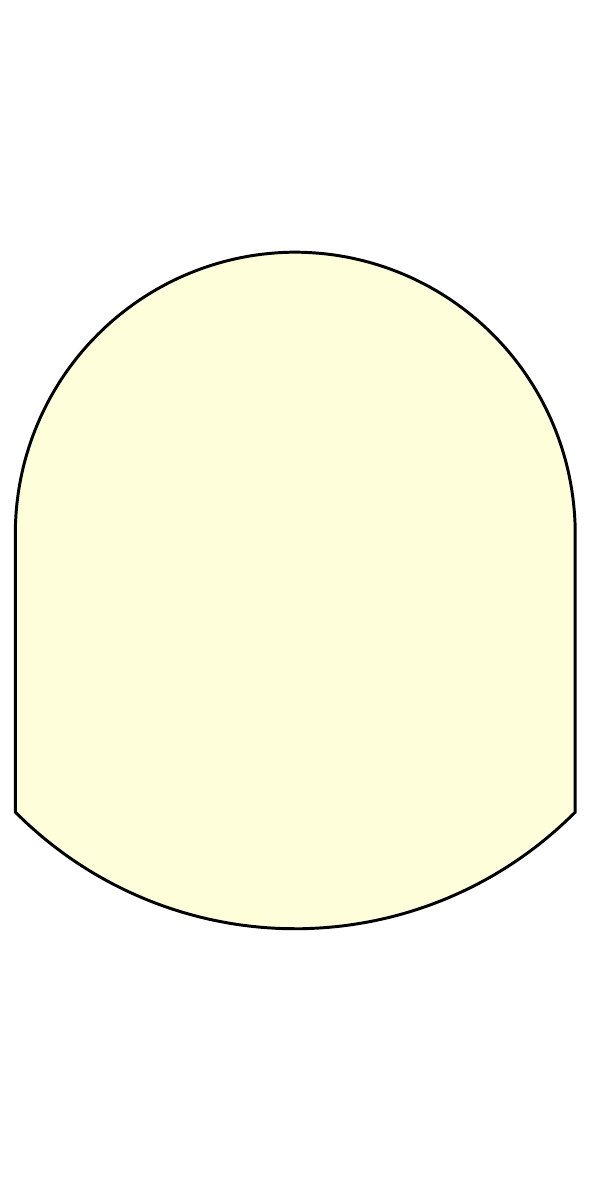}
\caption{$\psi=90^\circ$, $Area=1.0353...$}
\end{subfigure}

\caption{Representative shapes for Pattern~2}

\label{fig:second}

\end{figure}

To further illustrate the geometric characteristics of the optimal solutions, Figures~\ref{fig:first} and ~\ref{fig:second} show the locally optimal sofa shapes for several representative angles under the two motion patterns, respectively.

\clearpage

\clearpage

\section{Conclusion}

\begin{figure}[H]
\centering
\includegraphics[width=0.8\textwidth]{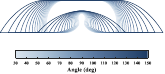}
\caption{Shape evolution with angle for Pattern~1}
\label{fig:jianbian1}
\end{figure}

\begin{figure}[H]
\centering
\includegraphics[width=0.8\textwidth]{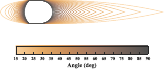}
\caption{Shape evolution with angle for Pattern~2}
\label{fig:jianbian2}
\end{figure}

We propose a gas-driven geometric optimization method and apply it to solve the maximal sofa problem for corridors with different angles. Through numerical simulations, we obtain locally optimal sofa shapes and areas under two motion patterns (see Figures~\ref{fig:jianbian1} and~\ref{fig:jianbian2}), and analyze their behavioral characteristics at various angles. In particular, the results reveal a phase transition: a critical point is identified at $43.327\ldots^\circ$, with a corresponding locally maximal area of $1.8674\ldots$, where the optimal motion pattern switches on either side of this point.

\noindent\textbf{Code and Data Availability}

The source code used for the numerical computations in this paper has been made publicly available. The code and the data required to reproduce the experiments can be obtained from the GitHub repository:

\begin{center}
\url{https://github.com/XingyiHe-rgb/moving-sofa-problem.git}
\end{center}

\section*{Acknowledgments}

The author thanks Professor Dan Romik for his invaluable guidance and encouragement.

\bibliographystyle{plain}
\bibliography{references}

@article{Hammersley1968,
  author  = {Hammersley, J. M.},
  title   = {On the enfeeblement of mathematical skills by ``Modern Mathematics'' and by similar soft intellectual trash in schools and universities},
  journal = {Bulletin of the Institute of Mathematics and Its Applications},
  volume  = {4},
  pages   = {66--85},
  year    = {1968}
}

@article{Gerver1992,
  author    = {Gerver, Joseph L.},
  title     = {On moving a sofa around a corner},
  journal   = {Geometriae Dedicata},
  volume    = {42},
  pages     = {267--283},
  year      = {1992},
  publisher = {Springer},
  doi       = {10.1007/BF02414066},
  url       = {https://doi.org/10.1007/BF02414066}
}

@article{Moser1966,
  author  = {Moser, Leo},
  title   = {Problem 66-11: Moving furniture through a hallway},
  journal = {SIAM Review},
  volume  = {8},
  number  = {3},
  pages   = {381},
  year    = {1966}
}

@misc{Gibbs2014,
  author       = {Gibbs, Philip},
  title        = {A Computational Study of Sofas and Cars},
  year         = {2014},
  note         = {Preprint},
  url          = {https://vixra.org/pdf/1411.0038v2.pdf}
}

@article{KallusRomik2018,
  author  = {Kallus, Yoav and Romik, Dan},
  title   = {Improved Upper Bounds in the Moving Sofa Problem},
  journal = {Advances in Mathematics},
  volume  = {340},
  pages   = {960--982},
  year    = {2018},
  doi     = {10.1016/j.aim.2018.10.022},
  url     = {https://arxiv.org/abs/1706.06630},
  eprint  = {1706.06630},
  archivePrefix = {arXiv},
  primaryClass  = {math}
}

@article{Romik2016,
  author  = {Romik, Dan},
  title   = {Differential Equations and Exact Solutions in the Moving Sofa Problem},
  journal = {Experimental Mathematics},
  volume  = {27},
  number  = {3},
  pages   = {299--313},
  year    = {2016},
  doi     = {10.1080/10586458.2016.1270858}
}

@article{Leng2024,
  author  = {Kuangdai Leng and Jia Bi and Jaehoon Cha and Samual Pingilla and Jeyan Thiyagalingam},
  title   = {Deep Learning Evidence for Global Optimality of Gerver's Sofa},
  journal = {Symmetry},
  volume  = {16},
  number  = {10},
  pages   = {1388},
  year    = {2024},
  doi     = {10.3390/sym16101388},
  url     = {https://www.mdpi.com/2073-8994/16/10/1388}
}

@misc{Deng2024,
  author       = {Deng, Zhipeng},
  title        = {Solving Moving Sofa Problem Using Calculus of Variations},
  year         = {2024},
  note         = {Preprint, arXiv:2407.02587 [math.CA]},
  url          = {http://arxiv.org/abs/2407.02587},
  eprint       = {2407.02587},
  archivePrefix= {arXiv},
  primaryClass = {math.CA}
}

@article{Maruyama1973,
  author  = {Maruyama, K.},
  title   = {An approximation method for solving the sofa problem},
  journal = {International Journal of Computer and Information Sciences},
  volume  = {2},
  pages   = {29--48},
  year    = {1973},
  doi     = {10.1007/BF00987154},
  url     = {https://doi.org/10.1007/BF00987154}
}

@article{Wagner1976,
  author  = {Wagner, N. R.},
  title   = {The sofa problem},
  journal = {American Mathematical Monthly},
  volume  = {83},
  number  = {3},
  pages   = {188--189},
  year    = {1976}
}

@misc{Georgiev2025,
  author       = {Georgiev, Bogdan and Gómez-Serrano, Javier and Tao, Terence and Wagner, Adam Zsolt},
  title        = {Mathematical Exploration and Discovery at Scale},
  year         = {2025},
  note         = {Preprint, arXiv:2511.02864 [cs.NE]},
  doi          = {10.48550/arXiv.2511.02864},
  url          = {https://arxiv.org/abs/2511.02864},
  eprint       = {2511.02864},
  archivePrefix= {arXiv},
  primaryClass = {cs.NE}
}

@misc{Baek2024,
  author       = {Baek, Jineon},
  title        = {Optimality of Gerver's Sofa},
  year         = {2024},
  note         = {Preprint, arXiv:2411.19826 [math.MG]},
  doi          = {10.48550/arXiv.2411.19826},
  url          = {https://arxiv.org/abs/2411.19826},
  eprint       = {2411.19826},
  archivePrefix= {arXiv},
  primaryClass = {math.MG}
}

@book{hadamard1923,
  title={Lectures on Cauchy's problem in linear partial differential equations},
  author={Hadamard, Jacques},
  year={1923},
  publisher={Yale University Press},
  address={New Haven},
  language={English}
}

\end{document}